\documentclass[11pt]{article}
\usepackage{latexsym}
\usepackage{amsfonts}
\usepackage{mathrsfs,amsthm,amsmath}
\usepackage{color}
\usepackage{todonotes}
\newtheorem{theorem}{Theorem}[section]
\newtheorem{proposition}[theorem]{Proposition}

\newtheorem{lemma}{Lemma}[section]

\newtheorem{assertion}{Assertion}[section]
\newtheorem{remark}{Remark}[section]
\newtheorem{example}{Example}[section]

\begin{document}
\title{Some examples of non-central moderate deviations for sequences of real random variables}
\author{Rita Giuliano\thanks{Address: Dipartimento di Matematica, Universit\`a di Pisa, 
Largo Bruno Pontecorvo 5, I-56127 Pisa, Italy. e-mail: \texttt{rita.giuliano@unipi.it}} 
\and Claudio Macci\thanks{Address: Dipartimento di Matematica, Universit\`a di Roma Tor 
Vergata, Via della Ricerca Scientifica, I-00133 Rome, Italy. e-mail:
\texttt{macci@mat.uniroma2.it}}}
\maketitle
\begin{abstract}
The term \emph{moderate deviations} is often used in the literature to mean a class of large 
deviation principles that, in some sense, fills the gap between a convergence in 
probability to zero (governed by a large deviation principle) and a weak convergence to a
centered Normal distribution. In this paper we present some examples of classes of large 
deviation principles of this kind, but the involved random variables converge weakly to 
Gumbel, exponential and Laplace distributions.
\end{abstract}

\section{Introduction}
The theory of large deviations gives an asymptotic computation of small probabilities on 
exponential scale; see \cite{DemboZeitouni} as a reference of this topic. The basic definition
of this theory is the concept of \emph{large deviation principle}, which provides some asymptotic 
bounds for a family of probability measures on the same topological space; these bounds are expressed 
in terms of a \emph{speed function} (that tends to infinity) and a lower semicontinuous \emph{rate 
function} defined on the topological space.

The term \emph{moderate deviations} is used for a class of large deviation principles which fill the
gap between a convergence to a constant (governed by a large deviation principle with a suitable speed 
function) and an asymptotic normality result. In view of the examples studied in this paper we explain 
the concept of moderate deviations in the next Assertion \ref{claim:MD}, and our presentation will be
restricted to sequences of real random variables defined on the same probability space 
$(\Omega,\mathcal{F},P)$; thus the speed function is a sequence $\{v_n:n\geq 1\}$ such that 
$v_n\to\infty$ (as in the rest of the paper).

\begin{assertion}[(Classical) moderate deviations]\label{claim:MD}
Let $\{C_n:n\geq 1\}$ be a sequence of real random variables such that the following asymptotic regimes
hold.\\
$\mathbf{R1}$: $\{C_n:n\geq 1\}$ converges in probability to zero, and this convergence is governed by
a large deviation principle with speed $v_n$ and rate function $I_{\mathrm{LD}}$ (such that 
$I_{\mathrm{LD}}(x)=0$ if and only if $x=0$);\\
$\mathbf{R2}$: $\sqrt{v_n}C_n$ converges weakly to a centered Normal distribution with (positive) 
variance $\sigma^2$.\\
Then we talk about moderate deviations when, for every family of positive numbers $\{a_n:n\geq 1\}$
such that
\begin{equation}\label{eq:MD-conditions}
a_n\to 0\ \mbox{and}\ a_nv_n\to\infty,
\end{equation}
the sequence of random variables $\{\sqrt{a_nv_n}C_n:n\geq 1\}$ satisfies the large deviation
principle with speed $1/a_n$ and rate function $I_{\mathrm{MD}}$ defined by
\begin{equation}\label{eq:MD-rf}
I_{\mathrm{MD}}(x)=\frac{x^2}{2\sigma^2}\quad\mbox{for all}\ x\in\mathbb{R}.
\end{equation}
Moreover one typically has 
\begin{equation}\label{eq:MD-rf-and-av}
I_{\mathrm{LD}}^{\prime\prime}(0)=\frac{1}{\sigma^2}.
\end{equation}
\end{assertion}

We have the following remarks.

\begin{remark}\label{rem:MD-fill-the-gap}
We can recover the asymptotic regimes $\mathbf{R1}$ and $\mathbf{R2}$ in Assertion \ref{claim:MD} by 
setting $a_n=\frac{1}{v_n}$ and $a_n=1$ respectively; note that, in both cases, one of the conditions
in \eqref{eq:MD-conditions} holds and the other one fails. So the class of large deviation principles 
in Assertion \ref{claim:MD} is determined by a family of positive scaling factors $\{a_n:n\geq 1\}$ 
that fill the gap between the asymptotic regimes $\mathbf{R1}$ and $\mathbf{R2}$. 
Moreover, by \eqref{eq:MD-conditions}, the speed $1/a_n$ for the random variables random variables 
$\{\sqrt{a_nv_n}C_n:n\geq 1\}$ has a lower intensity than the speed $v_n$ (see $\mathbf{R1}$).
\end{remark}

\begin{remark}\label{rem:auxiliary-centering}
Concerning the random variables $\{C_n:n\geq 1\}$ in Assertion \ref{claim:MD}, one often has
$$C_n=Z_n-z_\infty\quad\mbox{for all}\ n\geq 1,$$
where $\{Z_n:n\geq 1\}$ is a sequence which converges in probability to $z_\infty$ for some
$z_\infty\in\mathbb{R}$ (see for instance \eqref{eq:DD-paper-sequence} for Example \ref{ex:DD-paper}
where $z_\infty=t$); moreover this convergence is governed by a large deviation principle with some 
speed $v_n$, say, and rate function $I_Z$, and one has $I_Z(z)=0$ if and only if $z=z_\infty$. Then, 
for the rate function $I_{\mathrm{LD}}$ in Assertion \ref{claim:MD}, one has
$$I_{\mathrm{LD}}(x)=I_Z(x+z_\infty)\quad\mbox{for all}\ x\in\mathbb{R}.$$
Moreover, if we refer to the equality \eqref{eq:MD-rf-and-av} at the end of Assertion \ref{claim:MD},
one typically has $I_Z^{\prime\prime}(z_\infty)=\frac{1}{\sigma^2}$ because 
$I_Z^{\prime\prime}(z_\infty)=I_{\mathrm{LD}}^{\prime\prime}(0)$.
\end{remark}

Now we present a prototype example of the framework in Assertion \ref{claim:MD}, for which one
can refer to the law of large numbers and the central limit theorem. Without loss of generality we restrict
our attention to the case of centered random variables $\{X_n:n\geq 1\}$ in order to have a convergence
to zero (as in the asymptotic regime $\mathbf{R1}$ in Assertion \ref{claim:MD}).

\begin{example}\label{ex:prototype}
We set	
$$C_n:=\frac{X_1+\cdots+X_n}{n}\quad\mbox{for all}\ n\geq 1$$
where $\{X_n:n\geq 1\}$ is a sequence of i.i.d. centered real random variables. Moreover we also assume that 
$\mathbb{E}[e^{\theta X_1}]<\infty$ in a neighborhood of the origin $\theta=0$, and therefore
$\sigma^2=\mathrm{Var}[X_1]$ is finite (we also assume that $\sigma^2>0$ to avoid trivialities). In such
a case the sequence $\{C_n:n\geq 1\}$ satisfies the large deviation principle with speed $v_n=n$ and rate
function $I_{\mathrm{LD}}$ defined by
$$I_{\mathrm{LD}}(x)=\sup_{\theta\in\mathbb{R}}\left\{\theta x-\log\mathbb{E}\left[e^{\theta X_1}\right]\right\}\quad\mbox{for all}\ x\in\mathbb{R};$$
this is a consequence of Cram\'er Theorem (see e.g. Theorem 2.2.3 in \cite{DemboZeitouni}),
and one can check that $I_{\mathrm{LD}}^{\prime\prime}(0)=\frac{1}{\sigma^2}$. Moreover, for every family
of positive numbers $\{a_n:n\geq 1\}$ such that \eqref{eq:MD-conditions} holds, the sequence of random 
variables $\{\sqrt{a_nn}C_n:n\geq 1\}$ satisfies the large deviation principle with speed $1/a_n$ and 
rate function $I_{\mathrm{MD}}$ defined by \eqref{eq:MD-rf}; this is a consequence of Theorem 3.7.1 in 
\cite{DemboZeitouni} with $d=1$. Actually this prototype example can be presented for an arbitrary $d$ 
(i.e. for multivariate random variables) by taking into account the multivariate version of Cram\'er 
Theorem (see e.g. Theorem 2.2.30 in \cite{DemboZeitouni}).
\end{example}

The aim of this paper is to present some examples of \emph{non-central} moderate deviations. We use this 
terminology to mean a class of large deviation principles which fill the gap between a convergence to a 
constant (governed by a large deviation principle with some speed $v_n$ and some rate function which 
uniquely vanishes at that constant) and a weak convergence to some non-Gaussian law. This should happen 
in the same spirit of Assertion \ref{claim:MD} for some positive scalings $\{a_n:n\geq 1\}$ such that 
\eqref{eq:MD-conditions} holds. 

The terminology non-central moderate deviations appears in three recent results: Proposition 
3.3 in \cite{LeonenkoMacciPacchiarotti} (for possibly $m$-variate random variables), Proposition 4.3 in 
\cite{MacciMartinucciPirozzi} and Proposition 3.3 in \cite{BeghinMacci}. In the first two cases the weak 
convergence is trivial because one has a family of identically distributed random variables. Another result 
is Proposition 2.2 in \cite{IafrateMacci}, where the convergence in distribution is not trivial; however the 
term non-central moderate deviations does not appear in that paper.

The common line of the examples studied in this paper can be summarized as follows.

\begin{assertion}[Common line of several examples in this paper]\label{claim:ncMD-ex}
Let $\{C_n:n\geq 1\}$ be a sequence of real random variables such that the following asymptotic regimes hold.\\
$\mathbf{R1}$: $\{C_n:n\geq 1\}$ converges in probability to zero, and this convergence is governed by
a large deviation principle with speed $v_n$ (we always have $v_n\to\infty$) and rate function 
$I_{\mathrm{LD}}$ (such that $I_{\mathrm{LD}}(x)=0$ if and only if $x=0$);\\
$\mathbf{R2}$: $v_nC_n$ converges weakly to a non-Gaussian law.\\
Moreover, for every family of positive numbers $\{a_n:n\geq 1\}$ such that such that 
\eqref{eq:MD-conditions} holds, the sequence of random variables $\{a_nv_nC_n:n\geq 1\}$ satisfies
the large deviation principle with speed $1/a_n$ and a suitable rate function $I_{\mathrm{MD}}$.
\end{assertion}

Some interesting common features of the examples studied in this paper are given by the following
equalities: $I_{\mathrm{LD}}(0)=I_{\mathrm{MD}}(0)=0$ (as in Assertion \ref{claim:MD}),
$$I_{\mathrm{MD}}(x)=I_{\mathrm{LD}}(0+)x\ \mathrm{or}\ I_{\mathrm{MD}}(x)=\infty\quad\mbox{if}\ x>0$$
and
$$I_{\mathrm{MD}}(x)=I_{\mathrm{LD}}(0-)x\ \mathrm{or}\ I_{\mathrm{MD}}(x)=\infty\quad\mbox{if}\ x<0.$$
Not all the non-central moderate deviation results have these common features; for instance they do not
appear in Proposition 3.3 in \cite{BeghinMacci}. This explains the interest of non-central moderate
deviations that, in our opinion, deserve to be investigated. We also recall that the common features of
the examples studied in this paper can be seen as the analogue of the equality 
$I_{\mathrm{MD}}(x)=\frac{I_{\mathrm{LD}}^{\prime\prime}(0)}{2}x^2$ stated in Assertion \ref{claim:MD} 
for the classical moderate deviations (it is a consequence of \eqref{eq:MD-rf} and \eqref{eq:MD-rf-and-av}).

Now we present the outline of the paper with a very brief description of the examples studied in each
section. Section \ref{sec:preliminaries} is devoted to some preliminaries. In Section \ref{sec:minima}
we study an example with minima of i.i.d. nonnegative random variables. In Section \ref{sec:MDA-Gumbel} 
we study an example with maxima of i.i.d. random variables in the Maximum Domain Attraction of the Gumbel 
distribution (see e.g. the family of distributions in Theorem 8.13.4 in \cite{BinghamGoldieTeugels} 
together with the well-known Fisher-Tippett Theorem, e.g. Theorem 8.13.1 in \cite{BinghamGoldieTeugels});
another example with a weak convergence to the Gumbel distribution is studied in Section \ref{sec:occupancy},
and it concerns the classical occupancy problem (or the coupon collector's problem). Finally, in Section 
\ref{sec:DD-paper} we consider an example inspired by a recent replacement model for random lifetimes in the
literature (see e.g. \cite{DicrescenzoDigironimo}).

We conclude with some notation used throughout the paper. We write $a_n\sim b_n$ to mean that 
$\frac{a_n}{b_n}\to 1$ as $n\to\infty$; moreover we use the symbol $[x]$ for the integer part of 
$x\in\mathbb{R}$, i.e.
$$[x]:=\max\{k\in\mathbb{Z}:k\leq x\}.$$

\section{Preliminaries}\label{sec:preliminaries}
We start with the definition of large deviation principle (see e.g. \cite{DemboZeitouni}, pages
4-5). Let $(\mathcal{X},\tau_X)$ be a topological space and let $\{Z_n:n\geq 1\}$ be a sequence 
of $\mathcal{X}$-valued random variables defined on the same probability space $(\Omega,\mathcal{F},P)$.
A sequence $\{v_n:n\geq 1\}$ such that $v_n\to\infty$ (as $n\to\infty$) is called a \emph{speed function}, 
and a lower semicontinuous function $I:\mathcal{X}\to[0,\infty]$ is called a \emph{rate function}. Then the 
sequence $\{Z_n:n\geq 1\}$ satisfies the large deviation principle (LDP from now on) with speed $v_n$ and
rate function $I$ if both the following relations hold
\begin{equation}\label{eq:UB}
\limsup_{n\to\infty}\frac{1}{v_n}\log P(Z_n\in C)\leq-\inf_{x\in C}I(x)\quad\mbox{for all closed sets}\ C,
\end{equation}
and
\begin{equation}\label{eq:LB}
\liminf_{n\to\infty}\frac{1}{v_n}\log P(Z_n\in O)\geq-\inf_{x\in O}I(x)\quad\mbox{for all open sets}\ O.
\end{equation}
The rate function $I$ is said to be \emph{good} if, for every $\eta\geq 0$, the level set 
$\{x\in\mathcal{X}:I(x)\leq\eta\}$ is compact.

The following Lemma \ref{lemma} is quite a standard result (anyway we give briefly some hints of the proof
for the sake of completeness). All the moderate deviation results in this paper concern the situation 
presented in Assertion \ref{claim:ncMD-ex}, and they will be proved by applying Lemma \ref{lemma}; more precisely,
for every choice of the positive scalings $\{a_n:n\geq 1\}$ such that \eqref{eq:MD-conditions} holds (with a
further condition for Example \ref{ex:DD-paper}, i.e. \eqref{eq:restriction}), we shall have $s_n=1/a_n$ 
and $I=I_{\mathrm{MD}}$ for some rate function $I_{\mathrm{MD}}$. We shall apply Lemma \ref{lemma} also to 
prove Proposition \ref{prop:LD-DD-paper} (with $s_n=v_n=n$ and $I=I_{\mathrm{LD}}$ for a suitable rate function 
$I_{\mathrm{LD}}$), which provides the LDP for the asymptotic regime $\mathbf{R1}$ concerning Example 
\ref{ex:DD-paper}.

\begin{lemma}\label{lemma}
Let $\{s_n:n\geq 1\}$ be a speed and let $\{C_n:n\geq 1\}$ be a sequence of real random variables defined on 
the same probability space $(\Omega,\mathcal{F},P)$. Moreover let $I:\mathbb{R}\to[0,\infty]$ be a rate function
decreasing on $(-\infty,0)$, increasing on $(0,\infty)$ and such that $I(x)=0$ if and only if $x=0$. Moreover 
assume that:
\begin{equation}\label{eq:halfline1}
\limsup_{n\to\infty}\frac{1}{s_n}\log P(C_n\geq x)\leq-I(x)\quad\mbox{for all}\ x>0;
\end{equation}
\begin{equation}\label{eq:halfline2}
\limsup_{n\to\infty}\frac{1}{s_n}\log P(C_n\leq x)\leq-I(x)\quad\mbox{for all}\ x<0;
\end{equation}
\begin{equation}\label{eq:local-bound}
\liminf_{n\to\infty}\frac{1}{s_n}\log P(C_n\in O)\geq-I(x)\quad\left.\begin{array}{ll}
\mbox{for every $x\in\mathbb{R}$ such that $I(x)<\infty$}\\
\mbox{and for all open sets $O$ such that $x\in O$.}
\end{array}\right.
\end{equation}
Then $\{C_n:n\geq 1\}$ satisfies the LDP with speed $s_n$ and rate function $I$.
\end{lemma}
\begin{proof}
It is known that the final statement \eqref{eq:local-bound} yields the lower bound for open sets \eqref{eq:LB}
(see e.g. condition (b) with eq. (1.2.8) in \cite{DemboZeitouni}). Here we prove that \eqref{eq:halfline1} and
\eqref{eq:halfline2} yield the upper bound for closed sets \eqref{eq:UB}. Such an upper bound trivially holds 
if $C$ is the empty set or if $0\in C$, and therefore in what follows we assume that $C$ is nonempty and 
$0\notin C$; then at least one of the sets $C\cap (0,\infty)$ and $C\cap (-\infty,0)$ is nonempty. For simplicity
we assume that both sets $C\cap (0,\infty)$ and $C\cap (-\infty,0)$ are nonempty (in fact, if one of them is 
empty, the proof can be adapted readily). Then we can find $x_1\in C\cap (0,\infty)$ and 
$x_2\in C\cap (-\infty,0)$ such that $C\subset (-\infty,x_2]\cup[x_1,\infty)$, and we have
$$P(Z_n\in C)\leq P(Z_n\geq x_1)+P(Z_n\leq x_2);$$
thus, by Lemma 1.2.15 in \cite{DemboZeitouni} (together with \eqref{eq:halfline1} and \eqref{eq:halfline2}), we 
have
\begin{multline*}
\limsup_{n\to\infty}\frac{1}{s_n}\log P(Z_n\in C)
\leq\limsup_{n\to\infty}\frac{1}{s_n}\log\{P(Z_n\geq x_1)+P(Z_n\leq x_2)\}\\
=\max\{-I(x_1),-I(x_2)\}=-\min\{I(x_1),I(x_2)\}.
\end{multline*}
We conclude the proof noting that $\min\{I(x_1),I(x_2)\}=\inf_{x\in C}I(x)$ by the hypotheses.
\end{proof}

\section{An example with minima of i.i.d. nonnegative random variables}\label{sec:minima}
In this section we consider the following example.

\begin{example}\label{ex:minima}
Let $\{X_n:n\geq 1\}$ be a sequence of i.i.d. real random variables, with common distribution function $F$. We 
assume that $F$ is strictly increasing on $(\alpha_F,\omega_F)$, where
$$\alpha_F:=\inf\{x\in\mathbb{R}:F(x)>0\}=0\ \mbox{and}\ \omega_F:=\sup\{x\in\mathbb{R}:F(x)<1\};$$
so $\{X_n:n\geq 1\}$ are nonnegative random variables. We also assume that $F(0)=0$, and that there 
exists
\begin{equation}\label{eq:right-derivative-in-zero}
F^\prime(0+):=\lim_{x\to 0+}\frac{F(x)-F(0)}{x}=\lim_{x\to 0+}\frac{F(x)}{x}\in(0,\infty).
\end{equation}
In what follows we use the notation
$$F(x-):=\lim_{y\uparrow x}F(y).$$
\end{example}

Throughout this section we set
\begin{equation}\label{eq:minima-sequence}
C_n:=\min\{X_1,\ldots,X_n\}\quad\mbox{for all}\ n\geq 1.
\end{equation}

\begin{assertion}\label{claim:ncMD-ex-minima}
We can recover the asymptotic regimes $\mathbf{R1}$ and $\mathbf{R2}$ in Assertion \ref{claim:ncMD-ex} as
follows.\\
$\mathbf{R1}$: $\{C_n:n\geq 1\}$ converges in probability to zero (actually it is an a.s. convergence);
moreover, by Lemma 1 in \cite{MacciSMA}, $\{C_n:n\geq 1\}$ satisfies the LDP speed $v_n=n$ and rate 
function $I_{\mathrm{LD}}$ defined by
$$I_{\mathrm{LD}}(x):=\left\{\begin{array}{ll}
-\log(1-F(x-))&\ \mbox{if}\ x\in[0,\omega_F)\\
\infty&\ \mbox{otherwise}.
\end{array}\right.$$
$\mathbf{R2}$: $v_nC_n=n\min\{X_1,\ldots,X_n\}$ converges weakly to the exponential distribution with mean
$1/F^\prime(0+)$; in fact, for every $x>0$, by \eqref{eq:right-derivative-in-zero}, we get
\begin{multline*}
P(n\min\{X_1,\ldots,X_n\}\leq x)=1-P(n\min\{X_1,\ldots,X_n\}>x)=1-P^n\left(X_1>\frac{x}{n}\right)\\
=1-\left(1-F\left(\frac{x}{n}\right)\right)^n
=1-\left(1-\frac{nF\left(\frac{x}{n}\right)}{n}\right)^n\to 1-e^{-F^\prime(0+)x}\ \mbox{as}\ n\to\infty.
\end{multline*}
\end{assertion}

Now we prove the moderate deviation result.

\begin{proposition}\label{prop:MD-ex-minima}
For every family of positive numbers $\{a_n:n\geq 1\}$ such that \eqref{eq:MD-conditions} holds  
(i.e. $a_n\to 0$ and $a_nn\to\infty$), the sequence of random variables 
$\{a_nn\min\{X_1,\ldots,X_n\}:n\geq 1\}$ satisfies the LDP with speed $1/a_n$ and rate function 
$I_{\mathrm{MD}}$ defined by
$$I_{\mathrm{MD}}(x):=\left\{\begin{array}{ll}
F^\prime(0+)x&\ \mbox{if}\ x\in[0,\infty)\\
\infty&\ \mbox{otherwise}.
\end{array}\right.$$
\end{proposition}
\begin{proof}
We apply Lemma \ref{lemma} for every choice of the positive scalings $\{a_n:n\geq 1\}$ such that 
\eqref{eq:MD-conditions} holds, with $s_n=1/a_n$, $I=I_{\mathrm{MD}}$ and 
$\{C_n:n\geq 1\}$ as in \eqref{eq:minima-sequence}. 

\noindent\textsc{Proof of \eqref{eq:halfline1}.}  We have to show that, for every $x>0$,
$$\limsup_{n\to\infty}a_n\log P(a_nn\min\{X_1,\ldots,X_n\}\geq x)\leq-F^\prime(0+)x.$$
For every $\delta\in(0,x)$ we have
$$P(a_nn\min\{X_1,\ldots,X_n\}\geq x)=P^n\left(X_1\geq\frac{x}{a_nn}\right)=\left(1-F\left(\frac{x}{a_nn}-\right)\right)^n
\leq\left(1-F\left(\frac{x-\delta}{a_nn}\right)\right)^n$$
and, by the relation $\lim_{x\to 0+}F(x)=0$ and the limit in \eqref{eq:right-derivative-in-zero},
\begin{multline*}
\limsup_{n\to\infty}a_n\log P(a_nn\min\{X_1,\ldots,X_n\}\geq x)\\
\leq\limsup_{n\to\infty}a_nn\log\left(1-F\left(\frac{x-\delta}{a_nn}\right)\right)
=\limsup_{n\to\infty}a_nn\left(-F\left(\frac{x-\delta}{a_nn}\right)\right)
=-F^\prime(0+)(x-\delta);
\end{multline*}
so we obtain \eqref{eq:halfline1} by letting $\delta$ go to zero.

\noindent\textsc{Proof of \eqref{eq:halfline2}.} We have to show that, for every $x<0$,
$$\limsup_{n\to\infty}a_n\log P(a_nn\min\{X_1,\ldots,X_n\}\leq x)\leq-\infty.$$
This trivially holds as $-\infty\leq-\infty$ because the random variables $\{C_n:n\geq 1\}$ are 
nonnegative (and therefore $P(a_nn\min\{X_1,\ldots,X_n\}\leq x)=0$ for every $n\geq 1$).

\noindent\textsc{Proof of \eqref{eq:local-bound}.} We want to show that, for every $x\geq 0$ and 
for every open set $O$ such that $x\in O$, we have
$$\liminf_{n\to\infty}a_n\log P(a_nn\min\{X_1,\ldots,X_n\}\in O)\geq-F^\prime(0+)x.$$
The case $x=0$ is immediate; indeed, since $a_nn\min\{X_1,\ldots,X_n\}$ converges in probability to
zero by the Slutsky Theorem (by $a_n\to 0$ and the weak convergence in $\mathbf{R2}$ in Assertion 
\ref{claim:ncMD-ex-minima}), we have trivially $0\geq 0$.
So, from now on, we suppose that $x>0$ and we take $\delta>0$ small enough to have 
$(x-\delta,x+\delta)\subset O\cap(0,\infty)$. Then
\begin{multline*}
P(a_nn\min\{X_1,\ldots,X_n\}\in O)\geq P(a_nn\min\{X_1,\ldots,X_n\}\in(x-\delta,x+\delta))\\
=P\left(\min\{X_1,\ldots,X_n\}\in\left(\frac{x-\delta}{a_nn},\frac{x+\delta}{a_nn}\right)\right)
=P^n\left(X_1>\frac{x-\delta}{a_nn}\right)-P^n\left(X_1\geq\frac{x+\delta}{a_nn}\right)\\
=\left(1-F\left(\frac{x-\delta}{a_nn}\right)\right)^n-\left(1-F\left(\frac{x+\delta}{a_nn}-\right)\right)^n
=\left(1-F\left(\frac{x-\delta}{a_nn}\right)\right)^n
\left(1-\frac{\left(1-F\left(\frac{x+\delta}{a_nn}-\right)\right)^n}{\left(1-F\left(\frac{x-\delta}{a_nn}\right)\right)^n}\right)\\
=\left(1-F\left(\frac{x-\delta}{a_nn}\right)\right)^n
\left\{1-\exp\left(n\log\left(\frac{1-F\left(\frac{x+\delta}{a_nn}-\right)}{1-F\left(\frac{x-\delta}{a_nn}\right)}\right)\right)\right\}.
\end{multline*}
Now notice that
\begin{multline*}
\log\left(\frac{1-F\left(\frac{x+\delta}{a_nn}-\right)}{1-F\left(\frac{x-\delta}{a_nn}\right)}\right)
=\log\left(1+\frac{F\left(\frac{x-\delta}{a_nn}\right)-F\left(\frac{x+\delta}{a_nn}-\right)}{1-F\left(\frac{x-\delta}{a_nn}\right)}\right)\\
\leq\frac{F\left(\frac{x-\delta}{a_nn}\right)-F\left(\frac{x+\delta}{a_nn}-\right)}{1-F\left(\frac{x-\delta}{a_nn}\right)}
=-\frac{F\left(\frac{x+\delta}{a_nn}-\right)-F\left(\frac{x-\delta}{a_nn}\right)}{1-F\left(\frac{x-\delta}{a_nn}\right)}
\leq-\frac{F\left(\frac{x}{a_nn}\right)-F\left(\frac{x-\delta}{a_nn}\right)}{1-F\left(\frac{x-\delta}{a_nn}\right)}
\end{multline*}
and, again by the relation $\lim_{x\to 0+}F(x)=0$ and the limit in \eqref{eq:right-derivative-in-zero},
$$n\log\left(\frac{1-F\left(\frac{x+\delta}{a_nn}-\right)}{1-F\left(\frac{x-\delta}{a_nn}\right)}\right)
\leq-\frac{a_nnF\left(\frac{x}{a_nn}\right)-a_nnF\left(\frac{x-\delta}{a_nn}\right)}{a_n\left(1-F\left(\frac{x-\delta}{a_nn}\right)\right)}
\to-\infty\ \mbox{as}\ n\to\infty$$
because $a_n\left(1-F\left(\frac{x-\delta}{a_nn}\right)\right)\to 0$ and
$$a_nnF\left(\frac{x}{a_nn}\right)-a_nnF\left(\frac{x-\delta}{a_nn}\right)\to F^\prime(0+)(x-(x-\delta))=F^\prime(0+)\delta>0.$$
In conclusion we get
\begin{multline*}
\liminf_{n\to\infty}a_n\log	P(a_nn\min\{X_1,\ldots,X_n\}\in O)\\
\geq\liminf_{n\to\infty}a_n\log\left(1-F\left(\frac{x-\delta}{a_nn}\right)\right)^n
+\liminf_{n\to\infty}a_n\log\left(1-\exp\left(n\log\left(\frac{1-F\left(\frac{x+\delta}{a_nn}-\right)}
{1-F\left(\frac{x-\delta}{a_nn}\right)}\right)\right)\right)\\
=\liminf_{n\to\infty}a_nn\log\left(1-F\left(\frac{x-\delta}{a_nn}\right)\right)
=\liminf_{n\to\infty}a_nn\left(-F\left(\frac{x-\delta}{a_nn}\right)\right)=-F^\prime(0+)(x-\delta);
\end{multline*}
so we obtain \eqref{eq:local-bound} by letting $\delta$ go to zero.
\end{proof}

\section{An example with maxima of i.i.d. random variables in the MDA of Gumbel distribution}\label{sec:MDA-Gumbel}
In this section we consider the following example.

\begin{example}\label{ex:MDA-Gumbel}
Let $\{X_n:n\geq 1\}$ be a sequence of i.i.d. real random variables, with common distribution function $F$. 
Let $\omega_F$ be as in Example \ref{ex:minima}, and assume that $\omega_F=\infty$. We also assume that, 
for $x$ large enough, $F$ is strictly increasing with positive density $f$. Moreover let $w$ be the function
defined by
$$w(x):=\frac{\bar{F}(x)}{f(x)},\ \mbox{where}\ \bar{F}(x):=1-F(x),$$
and assume that $w$ is differentiable for $x$ large enough, $\lim_{x\to\infty}w'(x)=0$
and $w$ is a regularly varying function with exponent $1-\mu$ for some $\mu>0$, i.e.
$\lim_{x\to\infty}\frac{w(tx)}{w(x)}=t^{1-\mu}$ for all $t>0$. So it is well-known that
\begin{equation}\label{eq:w-representation-with-L}
w(x)=x^{1-\mu}L(x)
\end{equation}
for a suitable slowly varying function $L$, i.e. a function such that
$$\lim_{x\to\infty}\frac{L(tx)}{L(x)}=1\quad\mbox{for all}\ t>0$$
(this can be immediately checked by the definitions of regularly varying and slowly varying functions).
\end{example}

Here we list some particular cases in which $w$ is a regularly varying function with exponent $1-\mu$ for 
some $\mu>0$.

\begin{enumerate}
\item Standard normal distribution (see e.g. \cite{Small} (pp. 48, 50-51, 88-90) for the asymptotic behavior of $w(x)$)\\
\begin{tabular}{ccc}
$\bar{F}(x)=\frac{1}{\sqrt{2\pi}}\int_x^\infty e^{-t^2/2}dt$ ($x\in\mathbb{R}$),&\ $w(x)=\frac{\int_x^\infty e^{-t^2/2}dt}{e^{-x^2/2}}\sim\frac{1}{x}$,&\ $\mu=2$.
\end{tabular}
\item Gamma distribution for $a>0$ (see e.g. \cite{AbramowitzStegun} (formula 6.5.32, p. 263), for the limit of $w(x)$)\\
\begin{tabular}{ccc}
$\bar{F}(x)=\frac{1}{\Gamma(a)}\int_x^\infty t^{a-1}e^{-t}dt$ ($x>0$),&\ $w(x)=\frac{\int_x^\infty t^{a-1}e^{-t}dt}{x^{a-1}e^{-x}}\to 1$,&\ $\mu=1$.
\end{tabular}
\item Weibull distribution for $a>0$\\
\begin{tabular}{ccc}
$\bar{F}(x)=e^{-x^a}$ ($x>0$),&\ $w(x)=\frac{e^{-x^a}}{ax^{a-1}e^{-x^a}}=\frac{x^{1-a}}{a}$,&\ $\mu=a$.
\end{tabular}
\item Logistic distribution\\
\begin{tabular}{ccc}
$\bar{F}(x)=\frac{1}{1+e^x}$ ($x\in\mathbb{R}$),&\ $w(x)=\frac{\frac{1}{1+e^x}}{\frac{e^x}{(1+e^x)^2}}=1+e^{-x}\to 1$,&\ $\mu=1$.
\end{tabular}
\end{enumerate}

We need several consequences of the assumptions in Example \ref{ex:MDA-Gumbel}. We start with some results concerning the function
$L$ and the value $\mu$ (see \eqref{eq:w-representation-with-L}). In particular we provide an estimate for
\begin{equation}\label{eq:def-ell}
\ell:=\limsup_{x\to\infty}-\frac{L(x)\log\bar{F}(x)}{x^\mu}.
\end{equation}

\begin{lemma}\label{lem:consequence-of-Feller-result}
Under the assumptions in Example \ref{ex:MDA-Gumbel} we have $\lim_{x\to\infty}\frac{L(x)}{x^\mu}=0$.
\end{lemma}
\begin{proof}
It is known (see e.g. \cite{Feller}, Chapter VIII, Section 8, Lemma 2, p. 277) that,
for every $\varepsilon>0$, there exists $x_\varepsilon>0$ such that
$$x^{-\varepsilon}<L(x)<x^\varepsilon\ \mbox{for all}\ x>x_\varepsilon.$$
Then, if we take $\varepsilon<\mu$, we get
$$x^{-(\varepsilon+\mu)}<\frac{L(x)}{x^\mu}<x^{\varepsilon-\mu}\ \mbox{for all}\ x>x_\varepsilon,$$
and we immediately get the desired limit by letting $x$ go to infinity.
\end{proof}

\begin{lemma}\label{lem:Rita}
Under the assumptions in Example \ref{ex:MDA-Gumbel} we have $\ell\leq\frac{1}{\mu}$ (where $\ell$ is defined by 
\eqref{eq:def-ell}).
\end{lemma}
\begin{proof}
Recall from \eqref{eq:w-representation-with-L} that $\frac{\bar{F}(x)}{f(x)}=x^{1-\mu}L(x)$; thus $L(x)$ does not
vanish because we have $\bar{F}(x)\neq 0$ for all $x$ since $\omega_F=\infty$. Therefore, for some $x_1$, we get
$$\frac{f(x)}{\bar{F}(x)}=\frac{x^{\mu-1}}{L(x)}\ \mbox{for all}\ x>x_1,$$ 
which yields
$$\underbrace{\int_{x_1}^x\frac{f(t)}{\bar{F}(t)}dt}_{-\log\bar{F}(x)+\log\bar{F}(x_1)}=\int_{x_1}^x\frac{t^{\mu-1}}{L(t)}dt.$$
Recall also the so-called \emph{Potter bound} (see e.g. Theorem 1.5.6(i) in \cite{BinghamGoldieTeugels}): for every 
$A>1$ and $\delta>0$ there exists $x_0=x_0(A,\delta)$ such that
$$\frac{L(y)}{L(z)}\leq A\max\left\{\left(\frac{z}{y}\right)^\delta,\left(\frac{y}{z}\right)^\delta\right\}
\ \mbox{for all}\ y,z>x_0.$$
Then we have
$$-\frac{L(x)\log\bar{F}(x)}{x^\mu}+\frac{L(x)\log\bar{F}(x_1)}{x^\mu}
=\frac{L(x)}{x^\mu}\int_{x_1}^{x_0}\frac{t^{\mu-1}}{L(t)}dt
+\frac{L(x)}{x^\mu}\int_{x_0}^x\frac{t^{\mu-1}}{L(t)}dt;$$
moreover, by the Potter bound and some computations, we can estimate the last term as follows
$$\frac{L(x)}{x^\mu}\int_{x_0}^x\frac{t^{\mu-1}}{L(t)}dt
\leq\frac{A}{x^\mu}\int_{x_0}^x\left(\frac{x}{t}\right)^\delta t^{\mu-1}dt
=\frac{A}{\mu-\delta}\left(1-\left(\frac{x_0}{x}\right)^{\mu-\delta}\right).$$
Then, by the definition of $\ell$ in \eqref{eq:def-ell} and by Lemma \ref{lem:consequence-of-Feller-result}, we get 
$\ell\leq\frac{A}{\mu-\delta}$ by letting $x$ go to infinity; so we conclude by the arbitrariness of $A$ and $\delta$.
\end{proof}

\begin{remark}[A discussion on the inequality in Lemma \ref{lem:Rita}]\label{rem:Rita}
If we consider the four distributions listed above for which the function $w$ is regularly varying with exponent $1-\mu$ for 
some $\mu>0$, we have
\begin{equation}\label{eq:existence-of-L(infty)}
\lim_{x\to\infty}L(x)=L(\infty)\quad \mbox{for some}\ L(\infty)\in(0,\infty);
\end{equation}
indeed we have $L(\infty)=1$ for standard normal, Gamma and logistic distributions, and $L(\infty)=1/a$ for Weibull distribution.
Then, by the Hopital rule, we have
$$\lim_{x\to\infty}\frac{\log\bar{F}(x)}{x^\mu}=\lim_{x\to\infty}\frac{-f(x)/\bar{F}(x)}{\mu x^{\mu-1}}
=-\lim_{x\to\infty}\frac{1/w(x)}{\mu x^{\mu-1}}=-\lim_{x\to\infty}\frac{1}{\mu L(x)}=-\frac{1}{\mu L(\infty)},$$
and therefore
$$\lim_{x\to\infty}-\frac{L(x)\log\bar{F}(x)}{x^\mu}=-\lim_{x\to\infty}L(x)\lim_{x\to\infty}\frac{\log\bar{F}(x)}{x^\mu}=\frac{1}{\mu}.$$
In conclusion we have shown that the limit of $-\frac{L(x)\log\bar{F}(x)}{x^\mu}$ as $x\to\infty$ exists if 
\eqref{eq:existence-of-L(infty)} holds (as it happens for the four distributions listed above). We are not aware of cases in 
which the inequality in Lemma \ref{lem:Rita} is strict.
\end{remark}

Some further preliminaries are needed. Firstly, under the assumptions in Example \ref{ex:MDA-Gumbel}, the following
quantities are well-defined for $n$ large enough:
$$m_n:=F^{-1}\left(1-\frac{1}{n}\right)$$
and (in one of the following equalities we take into account \eqref{eq:w-representation-with-L})
\begin{equation}\label{eq:def-hn}
h_n:=m_nnf(m_n)=m_n\frac{f(m_n)}{\bar{F}(m_n)}=\frac{m_n}{w(m_n)}=\frac{m_n}{m_n^{1-\mu}L(m_n)}=\frac{m_n^\mu}{L(m_n)}.
\end{equation}
The following lemmas provide some properties of the function $w$ and the sequence $\{h_n:n\geq 1\}$.

\begin{lemma}\label{lem:MDA-Gumbel-w-property}
Under the assumptions in Example \ref{ex:MDA-Gumbel} we have $w(x_n)\sim w(y_n)$ for
$x_n,y_n\to\infty$ such that $x_n\sim y_n$.
\end{lemma}
\begin{proof}
By the well-known Karamata's representation of slowly varying functions (see e.g. Theorem 1.3.1 in
\cite{BinghamGoldieTeugels}), there exists $b>0$ such that the function $L$ introduced above can be 
written as
$$L(x)=c_1(x)\exp\left(\int_b^x\frac{c_2(t)}{t}dt\right)\quad\mbox{for all}\ x\geq b,$$
where $c_1$ and $c_2$ are suitable functions such that $c_1(x)$ tends to some finite limit and 
$c_2(x)\to 0$ (as $x\to\infty$). Then we have to prove that
$$\frac{w(x_n)}{w(y_n)}=\frac{c_1(x_n)\exp\left(\int_b^{x_n}\frac{c_2(t)}{t}dt\right)x_n^{1-\mu}}
{c_1(y_n)\exp\left(\int_b^{y_n}\frac{c_2(t)}{t}dt\right)y_n^{1-\mu}}\to 1\ \mbox{as}\ n\to\infty.$$
By the hypotheses we only need to prove that
$$\frac{\exp\left(\int_b^{x_n}\frac{c_2(t)}{t}dt\right)}{\exp\left(\int_b^{y_n}\frac{c_2(t)}{t}dt\right)}
=\exp\left(\int_{y_n}^{x_n}\frac{c_2(t)}{t}dt\right)\to 1,$$
which amounts to
$$\lim_{n\to\infty}\int_{y_n}^{x_n}\frac{c_2(t)}{t}dt=0.$$
Let $\delta\in(0,1)$ and $\varepsilon\in(0,\frac{\delta}{1+\delta})$ be arbitarily fixed.
Then there exists $n_\varepsilon\geq 1$ such that $|c_2(t)|<\varepsilon$ for $t>n_\varepsilon$ and
$1-\varepsilon<\frac{x_n}{y_n}<1+\varepsilon$ for $n>n_\varepsilon$. So, for $n$ large enough to have 
$x_n\wedge y_n>n_\varepsilon$, we have
$$\left|\int_{y_n}^{x_n}\frac{c_2(t)}{t}dt\right|\leq\int_{x_n\wedge y_n}^{x_n\vee y_n}\frac{|c_2(t)|}{t}dt
\leq\varepsilon\int_{x_n\wedge y_n}^{x_n\vee y_n}\frac{1}{t}dt=\varepsilon\log\frac{x_n\vee y_n}{x_n\wedge y_n}$$
and
$$1\leq\frac{x_n\vee y_n}{x_n\wedge y_n}=\left\{\begin{array}{ll}
\frac{x_n}{y_n}&\ \mbox{if}\ x_n\geq y_n\\
\frac{y_n}{x_n}&\ \mbox{if}\ x_n<y_n
\end{array}\right.\leq\left\{\begin{array}{ll}
1+\varepsilon&\ \mbox{if}\ x_n\geq y_n\\
\frac{1}{1-\varepsilon}&\ \mbox{if}\ x_n<y_n
\end{array}\right.<1+\delta.$$
The proof is complete.
\end{proof}

\begin{lemma}\label{lem:MDA-Gumbel-logarithmic-hn}
Under the assumptions in Example \ref{ex:MDA-Gumbel} we have $h_n\sim\mu\log n$ (as $n\to\infty$).
\end{lemma}
\begin{proof}
Firstly we recall that $h_n=\frac{m_n^\mu}{L(m_n)}$ (see \eqref{eq:def-hn}); then $h_n\to\infty$ by taking 
into account that $m_n\to\infty$ and by Lemma \ref{lem:consequence-of-Feller-result}. Moreover we have
$$\int_0^{m_n}\frac{1}{w(t)}dt=\int_0^{m_n}\frac{f(t)}{\bar{F}(t)}dt=-\int_0^{m_n}\frac{d}{dt}\log\bar{F}(t)dt
=-\log\bar{F}(m_n)+\log\bar{F}(0)=\log n+\log\bar{F}(0);$$
then $h_n\sim\mu\log n$ if and only if $h_n\sim\mu\int_0^{m_n}\frac{1}{w(t)}dt$, which amounts to
\begin{equation}\label{eq:key-limit}
\lim_{n\to\infty}\frac{\mu}{h_n}\int_0^{m_n}\frac{1}{w(t)}dt=1.
\end{equation}
So we prove the lemma showing that \eqref{eq:key-limit} holds. We take $\eta>0$ and we write
$$\frac{\mu}{h_n}\int_0^{m_n}\frac{1}{w(t)}dt=\underbrace{\frac{\mu}{h_n}\int_0^{\eta m_n}\frac{1}{w(t)}dt}_{=:I_n^{(1)}(\eta)}
+\underbrace{\frac{\mu}{h_n}\int_{\eta m_n}^{m_n}\frac{1}{w(t)}dt}_{=:I_n^{(2)}(\eta)}.$$
Now we estimate $I_n^{(1)}(\eta)$ and $I_n^{(2)}(\eta)$ separately.
\begin{itemize}
\item We have (make the change of variable $r=\bar{F}(t)$)
$$I_n^{(1)}(\eta)=\frac{\mu}{h_n}\int_0^{\eta m_n}\frac{f(t)}{\bar{F}(t)}dt
=\frac{\mu}{h_n}\int_{\bar{F}(\eta m_n)}^{\bar{F}(0)}\frac{dr}{r}=\frac{\mu}{h_n}(\log\bar{F}(0)-\log\bar{F}(\eta m_n)),$$
where $\bar{F}(0),\bar{F}(\eta m_n)<1$ because $\omega_F=\infty$; thus
$$0\leq\liminf_{n\to\infty}I_n^{(1)}(\eta)\leq\limsup_{n\to\infty}I_n^{(1)}(\eta)\leq\eta^\mu\mu\ell$$
because $\frac{\mu}{h_n}\log\bar{F}(0)\to 0$ (since $h_n\to\infty$), $-\frac{\mu}{h_n}\log\bar{F}(\eta m_n)\geq 0$,
\begin{multline*}
-\frac{\mu}{h_n}\log\bar{F}(\eta m_n)=-\frac{\mu}{m_n^\mu/L(m_n)}\log\bar{F}(\eta m_n)\\
=-\mu\eta^\mu\frac{L(m_n)}{(\eta m_n)^\mu}\log\bar{F}(\eta m_n)
=\mu\eta^\mu\frac{L(m_n)}{L(\eta m_n)}\left(-\frac{L(\eta m_n)\log\bar{F}(\eta m_n)}{(\eta m_n)^\mu}\right),
\end{multline*}
and by taking into account that $m_n\to\infty$, $L$ is a slowly varying function, and Lemma \ref{lem:Rita}.
\item We have (make the change of variable $u=\frac{t}{m_n}$)
$$I_n^{(2)}(\eta)=\frac{\mu}{m_n^\mu/L(m_n)}\int_{\eta m_n}^{m_n}\frac{1}{t^{1-\mu}L(t)}dt
=\frac{\mu L(m_n)}{m_n^\mu}\int_\eta^1\frac{m_ndu}{(m_nu)^{1-\mu}L(m_nu)}
=\int_\eta^1\frac{L(m_n)}{L(m_nu)}\mu u^{\mu-1}du,$$
and, since $\frac{L(m_n)}{L(m_nu)}\to 1$ uniformly on compact subsets of $(0,\infty)$ (with respect to $u$)
by Theorem 1.2.1 in \cite{BinghamGoldieTeugels} (here we again take into account that $m_n\to\infty$), for
all $\rho>0$ there exists $n_0=n_0(\eta,\rho)$ such that for all $n>n_0$
$$1-\rho<\frac{L(m_n)}{L(m_nu)}<1+\rho\quad\mbox{for all}\ u\in[\eta,1],$$
and therefore
$$(1-\rho)(1-\eta^\mu)<I_n^{(2)}(\eta)<(1+\rho)(1-\eta^\mu).$$
\end{itemize}
Finally we combine the estimates for $I_n^{(1)}(\eta)$ and $I_n^{(2)}(\eta)$, and we get \eqref{eq:key-limit}
by the arbitrariness of $\eta>0$ and $\rho>0$. This completes the proof.
\end{proof}

Throughout this section we set
\begin{equation}\label{eq:MDA-Gumbel-sequence}
C_n:=\frac{M_n}{m_n}-1\ \mbox{for all $n\geq n_0$, for some $n_0$, where $M_n:=\max\{X_1,\ldots,X_n\}$};
\end{equation}
we have to consider $n_0$ large enough in order to have a well-defined $m_n$ for $n\geq n_0$.

\begin{assertion}\label{claim:ncMD-ex-MDA-Gumbel}
We can recover the asymptotic regimes $\mathbf{R1}$ and $\mathbf{R2}$ in Assertion \ref{claim:ncMD-ex} as
follows.\\
$\mathbf{R1}$: $\{C_n:n\geq n_0\}$ converges in probability to zero (actually the authors have checked 
the almost sure convergence with an argument based on Theorem 4.4.4 in \cite{Galambos}, p. 268).
Moreover, as an immediate consequence of Proposition 3.1 in 
\cite{GiulianoMacciCSTM2014}, $\{C_n+1:n\geq n_0\}$ satisfies the LDP speed $v_n=\mu\log n$ (or 
equivalently $v_n=h_n$ by Lemma \ref{lem:MDA-Gumbel-logarithmic-hn}) and rate function $J$ defined by
$$J(y):=\left\{\begin{array}{ll}
\frac{y^\mu-1}{\mu}&\ \mbox{if}\ y\geq 1\\
\infty&\ \mbox{otherwise};
\end{array}\right.$$
then, since we deal with a sequence of shifted random variables, we deduce that 
$\{C_n:n\geq 2\}$ satisfies the LDP speed $v_n=h_n$ and rate function $I_{\mathrm{LD}}$ defined by
$$I_{\mathrm{LD}}(x):=J(x+1)=\left\{\begin{array}{ll}
\frac{(x+1)^\mu-1}{\mu}&\ \mbox{if}\ x\geq 0\\
\infty&\ \mbox{otherwise}.
\end{array}\right.$$
$\mathbf{R2}$: $v_nC_n=h_n\left(\frac{M_n}{m_n}-1\right)$ converges weakly to the Gumbel distribution by
a well-known result by von Mises (see e.g. Theorem 8.13.7 in \cite{BinghamGoldieTeugels}).
\end{assertion}

Now we prove the moderate deviation result. We shall see that, for this example, the rate functions
$I_{\mathrm{LD}}$ and $I_{\mathrm{MD}}$ coincide when $\mu=1$.

\begin{proposition}\label{prop:MD-ex-MDA-Gumbel}
For every family of positive numbers $\{a_n:n\geq 1\}$ such that \eqref{eq:MD-conditions} holds  
(i.e. $a_n\to 0$ and $a_nh_n\to\infty$), the sequence of random variables 
$\{a_nh_n(\frac{M_n}{m_n}-1):n\geq 1\}$ satisfies the LDP with speed $1/a_n$ and rate function
$I_{\mathrm{MD}}$ defined by
$$I_{\mathrm{MD}}(x):=\left\{\begin{array}{ll}
x&\ \mbox{if}\ x\in[0,\infty)\\
\infty&\ \mbox{otherwise}.
\end{array}\right.$$
\end{proposition}
\begin{proof}
We apply Lemma \ref{lemma} for every choice of the positive scalings $\{a_n:n\geq 1\}$ such that 
\eqref{eq:MD-conditions} holds, with $s_n=1/a_n$, $I=I_{\mathrm{MD}}$ and 
$\{C_n:n\geq n_0\}$ as in \eqref{eq:MDA-Gumbel-sequence}.

\noindent\textsc{Proof of \eqref{eq:halfline1}.} We have to show that, for every $x>0$,
$$\limsup_{n\to\infty}a_n\log P\left(a_nh_n\left(\frac{M_n}{m_n}-1\right)\geq x\right)\leq-x.$$
For $m_n^\prime:=(\frac{x}{a_nh_n}+1)m_n$ we get
\begin{multline*}
P\left(a_nh_n\left(\frac{M_n}{m_n}-1\right)\geq x\right)=P\left(M_n\geq m_n^\prime\right)
=1-P\left(M_n<m_n^\prime\right)=1-F^n\left(m_n^\prime\right)\\
=1-\left(1-\bar{F}\left(m_n^\prime\right)\right)^n
=1-\exp\left(n\log\left(1-\bar{F}\left(m_n^\prime\right)\right)\right)
\leq-n\log\left(1-\bar{F}\left(m_n^\prime\right)\right);
\end{multline*}
thus
\begin{multline*}
a_n\log	P\left(a_nh_n\left(\frac{M_n}{m_n}-1\right)\geq x\right)
\leq a_n\log\left(-n\log\left(1-\bar{F}\left(m_n^\prime\right)\right)\right)\\
=a_n\log n
+a_n\log\left(\frac{\log\left(1-\bar{F}\left(m_n^\prime\right)\right)}
{-\bar{F}\left(m_n^\prime\right)}\right)+a_n\log\bar{F}\left(m_n^\prime\right)
\end{multline*}
and, noting that
$$\lim_{n\to\infty}a_n\log\left(\frac{\log\left(1-\bar{F}\left(m_n^\prime\right)\right)}
{-\bar{F}\left(m_n^\prime\right)}\right)=0$$
because $m_n^\prime\to\infty$ (since $m_n\to\infty$), we obtain
$$\limsup_{n\to\infty}a_n\log P\left(a_nh_n\left(\frac{M_n}{m_n}-1\right)\geq x\right)\leq
\limsup_{n\to\infty}a_n\left(\log n+\log\bar{F}\left(m_n^\prime\right)\right).$$
Now we apply Lagrange Theorem (also known as the mean value theorem) to the function $g(z)=\log\bar{F}(z)$;
so there exists $\xi_n\in\left(m_n,m_n^\prime\right)$ such that
$$\log\bar{F}\left(m_n^\prime\right)=\log\bar{F}(m_n)-\frac{f(\xi_n)}{\bar{F}(\xi_n)}\left(m_n^\prime-m_n\right).$$
Then, by the definitions of $m_n$, $w$ and $h_n$, we get
\begin{multline*}
\log\bar{F}\left(m_n^\prime\right)=-\log n-\frac{m_nx}{w(\xi_n)a_nh_n}\\
=-\log n-\frac{x}{w(\xi_n)a_nnf(m_n)}
=-\log n-\frac{w(m_n)x}{w(\xi_n)a_nn\bar{F}(m_n)}=-\log n-\frac{w(m_n)x}{w(\xi_n)a_n},
\end{multline*}
and therefore
$$a_n\left(\log n+\log\bar{F}\left(m_n^\prime\right)\right)=-\frac{xw(m_n)}{w(\xi_n)}.$$
So we get
$$\limsup_{n\to\infty}a_n\log P\left(a_nh_n\left(\frac{M_n}{m_n}-1\right)\geq x\right)\leq
\limsup_{n\to\infty}a_n\left(\log n+\log\bar{F}\left(m_n^\prime\right)\right)\leq-x$$
by Lemma \ref{lem:MDA-Gumbel-w-property} with $x_n=m_n$ and $y_n=\xi_n$ (indeed
$\xi_n\sim m_n$ because $1\leq\frac{\xi_n}{m_n}\leq1+\frac{x}{a_nh_n}$ by construction, and
$1+\frac{x}{a_nh_n}\to 1$).
Thus \eqref{eq:halfline1} holds.

\noindent\textsc{Proof of \eqref{eq:halfline2}.} We have to show that, for every $x<0$,
$$\limsup_{n\to\infty}a_n\log P\left(a_nh_n\left(\frac{M_n}{m_n}-1\right)\leq x\right)\leq-\infty.$$
For $m_n^\prime:=(\frac{x}{a_nh_n}+1)m_n$ we get
$$P\left(a_nh_n\left(\frac{M_n}{m_n}-1\right)\leq x\right)=P\left(M_n\leq m_n^\prime\right)=F^n\left(m_n^\prime\right)
=\left(1-\bar{F}\left(m_n^\prime\right)\right)^n;$$
thus
\begin{multline*}
a_n\log	P\left(a_nh_n\left(\frac{M_n}{m_n}-1\right)\geq x\right)
\leq a_nn\log\left(1-\bar{F}\left(m_n^\prime\right)\right)\sim-a_nn\bar{F}\left(m_n^\prime\right)\\
=-\exp\left(\log\left(a_nn\bar{F}\left(m_n^\prime\right)\right)\right)
=-\exp\left(\log a_n+\log n+\log\bar{F}\left(m_n^\prime\right)\right)\\
=-\exp\left(\log a_n+\frac{a_n\left(\log n+\log\bar{F}\left(m_n^\prime\right)\right)}{a_n}\right).
\end{multline*}
Now we can repeat the computations above in the proof of \eqref{eq:halfline1} with some slight changes
(we mean the part with the application of Lagrange Theorem; the details are omitted), and we find
$$\lim_{n\to\infty}a_n\left(\log n+\log\bar{F}\left(m_n^\prime\right)\right)=-x>0,$$
and therefore
$$\lim_{n\to\infty}-\exp\left(\log a_n+\frac{a_n\left(\log n+\log\bar{F}\left(m_n^\prime\right)\right)}{a_n}\right)=-\infty.$$
Thus \eqref{eq:halfline2} holds.

\noindent\textsc{Proof of \eqref{eq:local-bound}.} We want to show that, for every $x\geq 0$ and 
for every open set $O$ such that $x\in O$, we have
$$\liminf_{n\to\infty}a_n\log P\left(a_nh_n\left(\frac{M_n}{m_n}-1\right)\in O\right)\geq-x.$$
The case $x=0$ is immediate; indeed, since $a_nh_n\left(\frac{M_n}{m_n}-1\right)$ converges in 
probability to zero by the Slutsky Theorem (by $a_n\to 0$ and the weak convergence in $\mathbf{R2}$
in Assertion \ref{claim:ncMD-ex-MDA-Gumbel}), we have trivially $0\geq 0$. So, from now on, we suppose
that $x>0$ and we take $\delta>0$ small enough to have $(x-\delta,x+\delta]\subset O\cap(0,\infty)$.
Then, for $m_n^{(\pm)}:=\left(1+\frac{x\pm\delta}{a_nh_n}\right)m_n$, we get
\begin{multline*}
P\left(a_nh_n\left(\frac{M_n}{m_n}-1\right)\in O\right)\geq	
P\left(x-\delta<a_nh_n\left(\frac{M_n}{m_n}-1\right)\leq x+\delta\right)\\
=P\left(m_n^{(-)}<M_n\leq m_n^{(+)}\right)=F^n\left(m_n^{(+)}\right)-F^n\left(m_n^{(-)}\right)\\
=\left(1-\bar{F}\left(m_n^{(+)}\right)\right)^n-\left(1-\bar{F}\left(m_n^{(-)}\right)\right)^n.
\end{multline*}
Now we apply Lagrange Theorem to the function $g(z)=(1-\bar{F}(z))^n$; so there exists
$\xi_n\in\left(m_n^{(-)},m_n^{(+)}\right)$
such that
$$\left(1-\bar{F}\left(m_n^{(+)}\right)\right)^n-\left(1-\bar{F}\left(m_n^{(-)}\right)\right)^n
=n\left(1-\bar{F}\left(\xi_n\right)\right)^{n-1}f(\xi_n)
\underbrace{\left(m_n^{(+)}-m_n^{(-)}\right)}_{=\frac{2\delta m_n}{a_nh_n}}.$$
Then, by the definition of $h_n$, we obtain
\begin{multline*}
a_n\log P\left(a_nh_n\left(\frac{M_n}{m_n}-1\right)\in O\right)\geq
a_n\left(\log n+(n-1)\log\left(1-\bar{F}\left(\xi_n\right)\right)+\log f(\xi_n)
+\log\frac{2\delta m_n}{a_nh_n}\right)\\
=a_n\left((n-1)\log\left(1-\bar{F}\left(\xi_n\right)\right)+\log f(\xi_n)
+\log(2\delta)-\log a_n-\log f(m_n)\right)
\end{multline*}
and therefore, by the definition of the function $w$,
\begin{multline*}
\liminf_{n\to\infty}a_n\log P\left(a_nh_n\left(\frac{M_n}{m_n}-1\right)\in O\right)
\geq\liminf_{n\to\infty}a_n\left((n-1)\log\left(1-\bar{F}\left(\xi_n\right)\right)+\log\frac{f(\xi_n)}{f(m_n)}\right)\\
=\liminf_{n\to\infty}a_n\left((n-1)\log\left(1-\bar{F}\left(\xi_n\right)\right)+\log\frac{w(m_n)\bar{F}(\xi_n)}{\bar{F}(m_n)w(\xi_n)}\right)\\
=\liminf_{n\to\infty}a_n\left((n-1)\log(1-\bar{F}(\xi_n))+\log\frac{w(m_n)}{w(\xi_n)}
+\log\frac{\bar{F}(\xi_n)}{\bar{F}(m_n)}\right);
\end{multline*}
then we get \eqref{eq:local-bound} if we show the following relations
\begin{equation}\label{eq:three-limits}
\lim_{n\to\infty}a_n\log\frac{w(m_n)}{w(\xi_n)}=0,\ \liminf_{n\to\infty}a_n\log\frac{\bar{F}(\xi_n)}{\bar{F}(m_n)}\geq-(x+\delta),
\ \lim_{n\to\infty}a_n(n-1)\log(1-\bar{F}(\xi_n))=0,
\end{equation}
and by letting $\delta$ go to zero. The first limit in \eqref{eq:three-limits} holds by Lemma 
\ref{lem:MDA-Gumbel-w-property} with $x_n=m_n$ and $y_n=\xi_n$ (indeed $\xi_n\sim m_n$
because $1+\frac{x-\delta}{a_nh_n}\leq\frac{\xi_n}{m_n}\leq1+\frac{x+\delta}{a_nh_n}$ 
by construction and $1+\frac{x\pm\delta}{a_nh_n}\to 1$). The inequality in \eqref{eq:three-limits}
holds since
$$a_n\log\frac{\bar{F}(\xi_n)}{\bar{F}(m_n)}\geq a_n\log\frac{\bar{F}(m_n^{(+)})}{\bar{F}(m_n)},$$
by a new application of Lagrange Theorem to the function $g(z)=-\log\bar{F}(z)$, and by 
applying again Lemma \ref{lem:MDA-Gumbel-w-property} (we omit the details to avoid repetitions).
Finally we prove the last limit in \eqref{eq:three-limits} if we show that
$$\lim_{n\to\infty}a_nn\bar{F}(m_n^{(\pm)})=0;$$
in fact
$$a_n(n-1)\log(1-\bar{F}(\xi_n))\sim -a_nn\bar{F}(\xi_n)\quad\mbox{and}\quad
a_nn\bar{F}(m_n^{(+)})\leq a_nn\bar{F}(\xi_n)\leq a_nn\bar{F}(m_n^{(-)}).$$
We have
$$a_nn\bar{F}(m_n^{(\pm)})=a_n\frac{\bar{F}(m_n^{(\pm)})}{\bar{F}(m_n)}
=\exp\left(\frac{a_n(\log\bar{F}(m_n^{(\pm)})-\log\bar{F}(m_n))}{a_n}+\log a_n\right);$$
moreover
$$\lim_{n\to\infty}a_n(\log\bar{F}(m_n^{(\pm)})-\log\bar{F}(m_n))=-(x\pm\delta)<0$$
(this follows once more from an application of Lagrange Theorem to the function $g(z)=-\log\bar{F}(z)$,
and by applying again Lemma \ref{lem:MDA-Gumbel-w-property}; the details are omitted). Thus
$$\lim_{n\to\infty}\frac{a_n(\log\bar{F}(m_n^{(\pm)})-\log\bar{F}(m_n))}{a_n}+\log a_n=-\infty,$$
which yields the desired last limit in \eqref{eq:three-limits}.
\end{proof}

\section{An example inspired by the classical occupancy problem}\label{sec:occupancy}
In this section we consider the following example.

\begin{example}\label{ex:occupancy}
Let $\{T_n:n\geq 1\}$ be a family of random variables such that
$$T_n:=\sum_{k=1}^nX_{n,k},$$
where $\{X_{n,k}:n\geq 1,1\leq k\leq n\}$ is a triangular array of random variables (i.e. 
$\{X_{n,k}:1\leq k\leq n\}$ are independent random variables for each fixed $n$), and each random
variable $X_{n,k}$ is geometric distributed with parameter $p_{n,k}:=1-\frac{k-1}{n}$, i.e.
$$P(X_{n,k}=h)=(1-p_{n,k})^{h-1}p_{n,k}\ \mbox{for all}\ h\geq 1.$$
\end{example}

It is well-known that the random variable $T_n$ can be seen as the number of balls required to fill $n$ 
boxes with at least one ball when one puts balls in $n$ boxes at random, and each ball is independently 
assigned to any fixed box with probability $\frac{1}{n}$; this is known in the literature as the classical 
occupancy problem. From a different point of view $T_n$ can also be related to the coupon collector's 
problem: a coupon collector chooses at random and independently among $n$ coupon types, and $T_n$ 
represents the number of coupons required to collect all the $n$ coupon types. Throughout this section 
we set
\begin{equation}\label{eq:occupancy-sequence}
C_n:=\frac{T_n}{n\log n}-1\ \mbox{for all}\ n\geq 2.
\end{equation}

\begin{assertion}\label{claim:ncMD-ex-occupancy}
We can recover the asymptotic regimes $\mathbf{R1}$ and $\mathbf{R2}$ in Assertion \ref{claim:ncMD-ex} as
follows.\\
$\mathbf{R1}$: $\{C_n:n\geq 2\}$ converges in probability to zero (see Example 2.2.7 in \cite{Durrett}
presented for the coupon collector's problem). Moreover, by Proposition 2.1 in \cite{GiulianoMacciTPMS2013},
$\{C_n+1:n\geq 2\}$ satisfies the LDP speed $v_n=\log n$ and rate function $J$ defined by
$$J(y):=\left\{\begin{array}{ll}
y-1&\ \mbox{if}\ y\geq 1\\
\infty&\ \mbox{otherwise};
\end{array}\right.$$
then, since we deal with a sequence of shifted random variables, we deduce that $\{C_n:n\geq 2\}$ satisfies 
the LDP speed $v_n=\log n$ and rate function $I_{\mathrm{LD}}$ defined by
$$I_{\mathrm{LD}}(x):=J(x+1)=\left\{\begin{array}{ll}
x&\ \mbox{if}\ x\geq 0\\
\infty&\ \mbox{otherwise}.
\end{array}\right.$$
$\mathbf{R2}$: $v_nC_n=\log n\left(\frac{T_n}{n\log n}-1\right)$ converges weakly to the Gumbel distribution 
(see e.g. Example 3.6.11 in \cite{Durrett}).
\end{assertion}

Now we prove the moderate deviation result. We shall see that, for this example, the rate functions
$I_{\mathrm{LD}}$ and $I_{\mathrm{MD}}$ coincide. Several parts of the proof of the next Proposition
\ref{prop:MD-ex-occupancy} have some analogies with the one presented for Proposition 2.1 in 
\cite{GiulianoMacciTPMS2013}; so we shall omit some details to avoid repetitions. 

\begin{proposition}\label{prop:MD-ex-occupancy}
For every family of positive numbers $\{a_n:n\geq 1\}$ such that \eqref{eq:MD-conditions} holds  
(i.e. $a_n\to 0$ and $a_n\log n\to\infty$), the sequence of random variables 
$\{a_n\log n(\frac{T_n}{n\log n}-1):n\geq 1\}$ satisfies the LDP with speed $1/a_n$ and rate function
$I_{\mathrm{MD}}$ defined by
$$I_{\mathrm{MD}}(x):=\left\{\begin{array}{ll}
x&\ \mbox{if}\ x\in[0,\infty)\\
\infty&\ \mbox{otherwise}.
\end{array}\right.$$
\end{proposition}
\begin{proof}
We apply Lemma \ref{lemma} for every choice of the positive scalings $\{a_n:n\geq 1\}$ such that 
\eqref{eq:MD-conditions} holds, with $s_n=1/a_n$, $I=I_{\mathrm{MD}}$ and 
$\{C_n:n\geq 2\}$ as in \eqref{eq:occupancy-sequence}.

\noindent\textsc{Proof of \eqref{eq:halfline1}.} We have to show that, for every $x>0$,
$$\limsup_{n\to\infty}a_n\log P\left(a_n\log n\left(\frac{T_n}{n\log n}-1\right)\geq x\right)\leq-x.$$
For every $\delta\in(0,x)$ we have (here the last inequality holds by a well-known estimate; see e.g.
Exercise 3.10 in \cite{MotwaniRaghavan}, page 58)
\begin{multline*}
P\left(a_n\log n\left(\frac{T_n}{n\log n}-1\right)\geq x\right)
=P\left(T_n\geq\left(\frac{x}{a_n\log n}+1\right)n\log n\right)\\
\leq P\left(T_n>\left(\frac{x-\delta}{a_n\log n}+1\right)n\log n\right)
\leq n^{1-(\frac{x-\delta}{a_n\log n}+1)}=n^{-\frac{x-\delta}{a_n\log n}}.
\end{multline*}
Thus
$$\limsup_{n\to\infty}a_n\log P\left(a_n\log n\left(\frac{T_n}{n\log n}-1\right)\geq x\right)
\leq\limsup_{n\to\infty}a_n\log n^{-\frac{x-\delta}{a_n\log n}}=-x+\delta,$$
and we obtain \eqref{eq:halfline1} by letting $\delta$ go to zero.

\noindent\textsc{Proof of \eqref{eq:halfline2}.} We have to show that, for every $x<0$,
$$\limsup_{n\to\infty}a_n\log P\left(a_n\log n\left(\frac{T_n}{n\log n}-1\right)\leq x\right)\leq-\infty.$$
We have $\frac{x}{a_n\log n}+1\in(0,1)$ eventually, and therefore (here the last inequality 
holds by a well-known estimate; see e.g. Theorem 5.10 and Corollary 5.11 in \cite{MitzenmacherUpfal})
\begin{multline*}
P\left(a_n\log n\left(\frac{T_n}{n\log n}-1\right)\leq x\right)=P\left(T_n\leq\left(\frac{x}{a_n\log n}+1\right)n\log n\right)\\
\leq 2\left(1-\exp\left(-\frac{\left(\frac{x}{a_n\log n}+1\right)n\log n}{n}\right)\right)^n
=2\left(1-n^{-\left(\frac{x}{a_n\log n}+1\right)}\right)^n.
\end{multline*}
Thus
\begin{multline*}
\limsup_{n\to\infty}a_n\log P\left(a_n\log n\left(\frac{T_n}{n\log n}-1\right)\leq x\right)
\leq\limsup_{n\to\infty}a_n\log\left(2\left(1-n^{-\left(\frac{x}{a_n\log n}+1\right)}\right)^n\right)\\
=\limsup_{n\to\infty}a_nn\log\left(1-n^{-\left(\frac{x}{a_n\log n}+1\right)}\right)
=\limsup_{n\to\infty}-a_nnn^{-\left(\frac{x}{a_n\log n}+1\right)}\\
=\limsup_{n\to\infty}-a_nn^{-\frac{x}{a_n\log n}}=\limsup_{n\to\infty}-a_ne^{-\frac{x}{a_n}}=-\infty,
\end{multline*}
and therefore \eqref{eq:halfline2} holds.

\noindent\textsc{Proof of \eqref{eq:local-bound}.} We want to show that, for every $x\geq 0$ and 
for every open set $O$ such that $x\in O$, we have
$$\liminf_{n\to\infty}a_n\log P\left(a_n\log n\left(\frac{T_n}{n\log n}-1\right)\in O\right)\geq-x.$$
The case $x=0$ is immediate; indeed, since $a_n\log n\left(\frac{T_n}{n\log n}-1\right)$ converges in 
probability to zero by the Slutsky Theorem (by $a_n\to 0$ and the weak convergence in $\mathbf{R2}$ in 
Assertion \ref{claim:ncMD-ex-occupancy}), we have trivially $0\geq 0$. So, from now on, we suppose that $x>0$ 
and we take $\delta>0$ small enough to have $(x-\delta,x+\delta]\subset O\cap(0,\infty)$. Moreover we also 
introduce the notation $F_{T_n}(\cdot)=P(T_n\leq\cdot)$ for the distribution function of $T_n$. Then
\begin{multline*}
P\left(a_n\log n\left(\frac{T_n}{n\log n}-1\right)\in O\right)\geq	
P\left(x-\delta<a_n\log n\left(\frac{T_n}{n\log n}-1\right)\leq x+\delta\right)\\
=P\left(\left(1+\frac{x-\delta}{a_n\log n}\right)n\log n<T_n\leq\left(1+\frac{x+\delta}{a_n\log n}\right)n\log n\right)\\
\geq F_{T_n}\left(\left[\left(1+\frac{x+\delta}{a_n\log n}\right)n\log n\right]\right)
-F_{T_n}\left(\left[\left(1+\frac{x-\delta}{a_n\log n}\right)n\log n\right]+1\right);
\end{multline*}
thus, by adapting some computations in the proof of Proposition 2.1 in \cite{GiulianoMacciTPMS2013},
we have
\begin{multline*}
F_{T_n}\left(\left[\left(1+\frac{x+\delta}{a_n\log n}\right)n\log n\right]\right)
-F_{T_n}\left(\left[\left(1+\frac{x-\delta}{a_n\log n}\right)n\log n\right]+1\right)\\
\geq A_n-(A_n^{(+)}+A_n^{(-)})=A_n\left(1-\frac{A_n^{(+)}+A_n^{(-)}}{A_n}\right)
\end{multline*}
where $A_n,A_n^{(+)},A_n^{(-)}$ are the following nonnegative quantities
\begin{equation}\label{eq:def-An}
A_n:=\left(1-e^{-\frac{1}{n}[(1+\frac{x+\delta}{a_n\log n})n\log n]}\right)^n
-\left(1-e^{-\frac{1}{n}([(1+\frac{x-\delta}{a_n\log n})n\log n]+1)}\right)^n,
\end{equation}
$$A_n^{(+)}:=\sum_{\mathrm{even}\ k}\binom{n}{k}e^{-\frac{k}{n}[(1+\frac{x+\delta}{a_n\log n})n\log n]}
\left(1-\left(1-\frac{k}{n}\right)^{\frac{k}{n}[(1+\frac{x+\delta}{a_n\log n})n\log n]}\right)$$
and
$$A_n^{(-)}:=\sum_{\mathrm{odd}\ k}\binom{n}{k}e^{-\frac{k}{n}([(1+\frac{x-\delta}{a_n\log n})n\log n]+1)}
\left(1-\left(1-\frac{k}{n}\right)^{\frac{k}{n}([(1+\frac{x-\delta}{a_n\log n})n\log n]+1)}\right).$$
Thus we can say that
$$a_n\log P\left(a_n\log n\left(\frac{T_n}{n\log n}-1\right)\in O\right)
\geq a_n\log A_n+a_n\log\left(1-\frac{A_n^{(+)}+A_n^{(-)}}{A_n}\right),$$
and therefore
$$\liminf_{n\to\infty}a_n\log P\left(a_n\log n\left(\frac{T_n}{n\log n}-1\right)\in O\right)
\geq\liminf_{n\to\infty}a_n\log A_n+\liminf_{n\to\infty}a_n\log\left(1-\frac{A_n^{(+)}+A_n^{(-)}}{A_n}\right).$$
We conclude the proof of \eqref{eq:local-bound} (and therefore the proof of the proposition)
if we show that
\begin{equation}\label{eq:condition(i)}
\liminf_{n\to\infty}a_n\log A_n\geq-(x+\delta)
\end{equation}
and
\begin{equation}\label{eq:conditions(ii)-(iii)}
\lim_{n\to\infty}\frac{A_n^{(+)}}{A_n}=0,\quad\lim_{n\to\infty}\frac{A_n^{(-)}}{A_n}=0;
\end{equation}
in fact we obtain \eqref{eq:local-bound} by letting $\delta$ go to zero.

\noindent\textsc{Proof of \eqref{eq:condition(i)}.} Concerning $A_n$ in \eqref{eq:def-An} we 
apply Lagrange Theorem to the function $g(z)=(1-e^{-z})^n$; so there exists
$\xi_n\in\left(\frac{([(1+\frac{x-\delta}{a_n\log n})n\log n]+1)}{n},\frac{[(1+\frac{x+\delta}{a_n\log n})n\log n]}{n}\right)$
such that
$$A_n=n(1-e^{-\xi_n})^{n-1}e^{-\xi_n}\left(\frac{[(1+\frac{x+\delta}{a_n\log n})n\log n]}{n}-\frac{([(1+\frac{x-\delta}{a_n\log n})n\log n]+1)}{n}\right).$$
Moreover
\begin{multline*}
(1-e^{-\xi_n})^{n-1}e^{-\xi_n}\geq\left(1-e^{-\frac{([(1+\frac{x-\delta}{a_n\log n})n\log n]+1)}{n}}\right)^{n-1}
e^{-\frac{[(1+\frac{x+\delta}{a_n\log n})n\log n]}{n}}\\
\geq\left(1-e^{-(1+\frac{x-\delta}{a_n\log n})\log n}\right)^{n-1}e^{-(1+\frac{x+\delta}{a_n\log n})\log n}
=\left(1-\frac{e^{-\frac{x-\delta}{a_n}}}{n}\right)^{n-1}\frac{e^{-\frac{x+\delta}{a_n}}}{n}
\end{multline*}
and
\begin{multline*}
\frac{[(1+\frac{x+\delta}{a_n\log n})n\log n]}{n}-\frac{([(1+\frac{x-\delta}{a_n\log n})n\log n]+1)}{n}
=\frac{[(1+\frac{x+\delta}{a_n\log n})n\log n]-1-[(1+\frac{x-\delta}{a_n\log n})n\log n]}{n}\\
\geq\frac{(1+\frac{x+\delta}{a_n\log n})n\log n-1-1-(1+\frac{x-\delta}{a_n\log n})n\log n}{n}
=-\frac{2}{n}+\frac{2\delta}{a_n};
\end{multline*}
then
\begin{multline*}
a_n\log A_n\geq a_n\log\left(\left(1-\frac{e^{-\frac{x-\delta}{a_n}}}{n}\right)^{n-1}e^{-\frac{x+\delta}{a_n}}\left(-\frac{2}{n}+\frac{2\delta}{a_n}\right)\right)\\
=a_n(n-1)\log\left(1-\frac{e^{-\frac{x-\delta}{a_n}}}{n}\right)-(x+\delta)+a_n\log\left(-\frac{2}{n}+\frac{2\delta}{a_n}\right)
\end{multline*}
whence we obtain \eqref{eq:condition(i)} from
$$\lim_{n\to\infty}a_n\log\left(-\frac{2}{n}+\frac{2\delta}{a_n}\right)=0$$
and
$$a_n(n-1)\log\left(1-\frac{e^{-\frac{x-\delta}{a_n}}}{n}\right)\sim -a_n(n-1)\frac{e^{-\frac{x-\delta}{a_n}}}{n}.
\sim-a_ne^{-\frac{x-\delta}{a_n}}\to 0.$$

\noindent\textsc{Proof of \eqref{eq:conditions(ii)-(iii)}.} We remark that
$$\lim_{n\to\infty}\left(1-e^{-\frac{1}{n}[(1+\frac{x+\delta}{a_n\log n})n\log n]}\right)^n=1
\quad\mbox{and}\quad\lim_{n\to\infty}\left(1-e^{-\frac{1}{n}([(1+\frac{x-\delta}{a_n\log n})n\log n]+1)}\right)^n=1$$
because
\begin{multline*}
1=\liminf_{n\to\infty}\left(1-\frac{e^{-\frac{x+\delta}{a_n}+\frac{1}{n}}}{n}\right)^n
\leq\liminf_{n\to\infty}\left(1-e^{-\frac{1}{n}[(1+\frac{x+\delta}{a_n\log n})n\log n]}\right)^n\\
\leq\limsup_{n\to\infty}\left(1-e^{-\frac{1}{n}[(1+\frac{x+\delta}{a_n\log n})n\log n]}\right)^n
\leq\limsup_{n\to\infty}\left(1-\frac{e^{-\frac{x+\delta}{a_n}}}{n}\right)^n=1
\end{multline*}
and the second limit can be proved with a similar argument.
Therefore
$$A_n=\left(1-e^{-\frac{1}{n}[(1+\frac{x+\delta}{a_n\log n})n\log n]}\right)^n
-\left(1-e^{-\frac{1}{n}([(1+\frac{x-\delta}{a_n\log n})n\log n]+1)}\right)^n\to 1-1=0.$$
Moreover
\begin{multline*}
A_n=\left(1-e^{-\frac{1}{n}([(1+\frac{x-\delta}{a_n\log n})n\log n]+1)}\right)^n
\left(\frac{\left(1-e^{-\frac{1}{n}[(1+\frac{x+\delta}{a_n\log n})n\log n]}\right)^n}
{\left(1-e^{-\frac{1}{n}([(1+\frac{x-\delta}{a_n\log n})n\log n]+1)}\right)^n}-1\right)\\
\sim\frac{\left(1-e^{-\frac{1}{n}[(1+\frac{x+\delta}{a_n\log n})n\log n]}\right)^n}
{\left(1-e^{-\frac{1}{n}([(1+\frac{x-\delta}{a_n\log n})n\log n]+1)}\right)^n}-1
=\exp\left(n\log\left(\frac{1-e^{-\frac{1}{n}[(1+\frac{x+\delta}{a_n\log n})n\log n]}}
{1-e^{-\frac{1}{n}([(1+\frac{x-\delta}{a_n\log n})n\log n]+1)}}\right)\right)-1,
\end{multline*}
where
\begin{equation}\label{eq:useful-limit}
\lim_{n\to\infty}n\log\left(\frac{1-e^{-\frac{1}{n}[(1+\frac{x+\delta}{a_n\log n})n\log n]}}
{1-e^{-\frac{1}{n}([(1+\frac{x-\delta}{a_n\log n})n\log n]+1)}}\right)=0.
\end{equation}
(the limit \eqref{eq:useful-limit} can be checked with some tedious computations:
$$n\log\left(\frac{1-e^{-\frac{1}{n}[(1+\frac{x+\delta}{a_n\log n})n\log n]}}
{1-e^{-\frac{1}{n}([(1+\frac{x-\delta}{a_n\log n})n\log n]+1)}}\right)
=n\log\left(1+\frac{e^{-\frac{1}{n}([(1+\frac{x-\delta}{a_n\log n})n\log n]+1)}
-e^{-\frac{1}{n}[(1+\frac{x+\delta}{a_n\log n})n\log n]}}
{1-e^{-\frac{1}{n}([(1+\frac{x-\delta}{a_n\log n})n\log n]+1)}}\right)$$
where 
$$e^{-\frac{1}{n}([(1+\frac{x-\delta}{a_n\log n})n\log n]+1)},e^{-\frac{1}{n}[(1+\frac{x+\delta}{a_n\log n})n\log n]}\to 0,$$
and therefore
\begin{multline*}
n\log\left(\frac{1-e^{-\frac{1}{n}[(1+\frac{x+\delta}{a_n\log n})n\log n]}}
{1-e^{-\frac{1}{n}([(1+\frac{x-\delta}{a_n\log n})n\log n]+1)}}\right)\sim
n\left(e^{-\frac{1}{n}([(1+\frac{x-\delta}{a_n\log n})n\log n]+1)}
-e^{-\frac{1}{n}[(1+\frac{x+\delta}{a_n\log n})n\log n]}\right)\\
=e^{\log n-\frac{1}{n}([(1+\frac{x-\delta}{a_n\log n})n\log n]+1)}
\left(1-e^{-\frac{1}{n}[(1+\frac{x+\delta}{a_n\log n})n\log n]
+\frac{1}{n}([(1+\frac{x-\delta}{a_n\log n})n\log n]+1)}\right);
\end{multline*}
moreover
$$e^{\log n-\frac{1}{n}([(1+\frac{x-\delta}{a_n\log n})n\log n]+1)},
e^{-\frac{1}{n}[(1+\frac{x+\delta}{a_n\log n})n\log n]
+\frac{1}{n}([(1+\frac{x-\delta}{a_n\log n})n\log n]+1)}\to 0
\quad\mbox{as}\ n\to\infty,$$
whence we obtain \eqref{eq:useful-limit}). Now, coming back to $A_n$, by 
\eqref{eq:useful-limit} we can say that
$$A_n\sim n\log\left(\frac{1-e^{-\frac{1}{n}[(1+\frac{x+\delta}{a_n\log n})n\log n]}}
{1-e^{-\frac{1}{n}([(1+\frac{x-\delta}{a_n\log n})n\log n]+1)}}\right);$$
moreover, by using some manipulations above for 
$\frac{1-e^{-\frac{1}{n}[(1+\frac{x+\delta}{a_n\log n})n\log n]}}
{1-e^{-\frac{1}{n}([(1+\frac{x-\delta}{a_n\log n})n\log n]+1)}}$, we get
\begin{multline*}
A_n\sim n\left(e^{-\frac{1}{n}([(1+\frac{x-\delta}{a_n\log n})n\log n]+1)}
-e^{-\frac{1}{n}[(1+\frac{x+\delta}{a_n\log n})n\log n]}\right)\\
=ne^{-\frac{1}{n}[(1+\frac{x+\delta}{a_n\log n})n\log n]}
\left(e^{-\frac{1}{n}([(1+\frac{x-\delta}{a_n\log n})n\log n]+1-[(1+\frac{x+\delta}{a_n\log n})n\log n])}-1\right).
\end{multline*}
Then
\begin{multline*}
A_n\sim ne^{-\frac{1}{n}[(1+\frac{x+\delta}{a_n\log n})n\log n]}
\left(e^{-\frac{1}{n}([(1+\frac{x-\delta}{a_n\log n})n\log n]+1-[(1+\frac{x+\delta}{a_n\log n})n\log n])}-1\right)\\
\geq ne^{-\frac{1}{n}(1+\frac{x+\delta}{a_n\log n})n\log n}
\left(e^{-\frac{1}{n}((1+\frac{x-\delta}{a_n\log n})n\log n+1-((1+\frac{x+\delta}{a_n\log n})n\log n-1))}-1\right)\\
=n\frac{1}{n}e^{-\frac{x+\delta}{a_n}}\left(e^{-\frac{2}{n}+\frac{2\delta}{a_n}}-1\right)
=e^{-\frac{x+\delta}{a_n}}\left(e^{-\frac{2}{n}+\frac{2\delta}{a_n}}-1\right)=:\tilde{A}_n.
\end{multline*}

Now, by using the same computations as in the proof of Proposition 2.1 in \cite{GiulianoMacciTPMS2013}, there exists 
a constant $C>0$ such that
$$0\leq A_n^{(+)}\leq\frac{C}{n}\left[\left(1+\frac{x+\delta}{a_n\log n}\right)n\log n\right]e^{-\frac{1}{n}[(1+\frac{x+\delta}{a_n\log n})n\log n]}
(1+e^{-\frac{1}{n}[(1+\frac{x+\delta}{a_n\log n})n\log n]})^{n-2}$$
and
$$0\leq A_n^{(-)}\leq\frac{C}{n}\left(\left[\left(1+\frac{x-\delta}{a_n\log n}\right)n\log n\right]+1\right)
e^{-\frac{1}{n}([(1+\frac{x-\delta}{a_n\log n})n\log n]+1)}
(1+e^{-\frac{1}{n}([(1+\frac{x-\delta}{a_n\log n})n\log n]+1)})^{n-2}.$$
Then, by observing $a_nn=a_n\log n\frac{n}{\log n}\to\infty$ and that
$$\lim_{n\to\infty}(1+e^{-\frac{1}{n}[(1+\frac{x+\delta}{a_n\log n})n\log n]})^{n-2}
=\lim_{n\to\infty}(1+e^{-\frac{1}{n}([(1+\frac{x-\delta}{a_n\log n})n\log n]+1)})^{n-2}=1,$$
we get
$$0\leq\limsup_{n\to\infty}\frac{A_n^{(+)}}{A_n}
\leq\limsup_{n\to\infty}\frac{C(\log n+\frac{x+\delta}{a_n})\frac{1}{n}e^{-\frac{x+\delta}{a_n}}}{\tilde{A}_n}
=\limsup_{n\to\infty}\frac{C(\log n+\frac{x+\delta}{a_n})\frac{1}{n}}{e^{-\frac{2}{n}+\frac{2\delta}{a_n}}-1}=0$$
and
$$0\leq\limsup_{n\to\infty}\frac{A_n^{(-)}}{A_n}
\leq\limsup_{n\to\infty}\frac{C(\log n+\frac{x-\delta}{a_n}+\frac{1}{n})\frac{1}{n}e^{-\frac{x-\delta}{a_n}}}{\tilde{A}_n}
=\limsup_{n\to\infty}\frac{C(\log n+\frac{x-\delta}{a_n}+\frac{1}{n})\frac{1}{n}}{e^{-\frac{2\delta}{a_n}}(e^{-\frac{2}{n}+\frac{2\delta}{a_n}}-1)}=0.$$
Thus the limits in \eqref{eq:conditions(ii)-(iii)} are checked.
\end{proof}

\section{An example inspired by the replacement model in \cite{DicrescenzoDigironimo}}\label{sec:DD-paper}
In this section we consider the following example.

\begin{example}\label{ex:DD-paper}
Let $F$ and $G$ be two continuous distribution functions on $\mathbb{R}$ such that $F(0)=G(0)=0$,
and assume that they are strictly increasing on $[0,\infty)$. Moreover we assume that, for some $t>0$,
there exist $F^\prime(t-)$ and $G^\prime(t+)$, i.e. the left derivative of $F(x)$ at $x=t$ and the right 
derivative of $G(x)$ at $x=t$, and that $F^\prime(t-),G^\prime(t+)>0$. Then let $\{Z_n:n\geq 1\}$ 
be a family of random variables defined on the same probability space $(\Omega,\mathcal{F},P)$ 
such that, for some $t>0$ and $\beta\in(0,1)$, their distribution functions are defined by
$$P(Z_n\leq x):=\left\{\begin{array}{ll}
\beta\left(\frac{F(x)}{F(t)}\right)^n&\ \mbox{if}\ x\in [0,t]\\
1-(1-\beta)\left(\frac{1-G(x)}{1-G(t)}\right)^n&\ \mbox{if}\ x\in (t,\infty).
\end{array}\right.$$
\end{example}

Note that, for every $n\geq 1$, the distribution function $P(Z_n\leq\cdot)$ is continuous. Moreover, if 
$\beta=F(t)$, after some easy computations one can check that
$$P(Z_1\leq x):=\left\{\begin{array}{ll}
F(x)&\ \mbox{if}\ x\in [0,t]\\
F(t)+\frac{1-F(t)}{1-G(t)}(G(x)-G(t))&\ \mbox{if}\ x\in (t,\infty);
\end{array}\right.$$
then $Z_1$ is the random lifetime that appears in a replacement model recently studied in the literature (see eq. (5) in 
\cite{DicrescenzoDigironimo}) where, at a fixed time $t$, an item with lifetime distribution function $F$ is
replaced by another item having the same age and a lifetime distribution function $G$. 

In general the distribution functions of the random variables $\{Z_n:n\geq 1\}$ are suitable mo\-difications of the 
distribution function of $Z_1$ and, as $n$ increases, the distribution of $Z_n$ is more concentrated around the point $t$. 
Throughout this section we set
\begin{equation}\label{eq:DD-paper-sequence}
C_n:=Z_n-t\ \mbox{for all}\ n\geq 1.
\end{equation}

\begin{assertion}\label{claim:ncMD-ex-DD-paper}
We can recover the asymptotic regimes $\mathbf{R1}$ and $\mathbf{R2}$ in Assertion \ref{claim:ncMD-ex} as
follows.\\
$\mathbf{R1}$: $\{C_n:n\geq 1\}$ converges in probability to zero because $\{Z_n:n\geq 1\}$ converges in probability
to $t$; in fact we have
$$\lim_{n\to\infty}P(Z_n\leq x)=\left\{\begin{array}{ll}
0&\ \mbox{if}\ x\in [0,t)\\
\beta&\ \mbox{if}\ x=t\\
1&\ \mbox{if}\ x\in (t,\infty).
\end{array}\right.$$
Moreover $\{C_n:n\geq 1\}$ satisfies a suitable LDP, with speed $v_n=n$, that will be presented in Proposition 
\ref{prop:LD-DD-paper} below.\\
$\mathbf{R2}$: $v_nC_n=n(Z_n-t)$ converges weakly to a suitable asymmetric Laplace distribution with distribution function
$H$ (see below). In fact we have
\begin{multline*}
P(n(Z_n-t)\leq x)=P\left(Z_n\leq t+\frac{x}{n}\right)
=\left\{\begin{array}{ll}
\beta\left(\frac{F(t+\frac{x}{n})}{F(t)}\right)^n&\ \mbox{if}\ x\leq 0\\
1-(1-\beta)\left(\frac{1-G(t+\frac{x}{n})}{1-G(t)}\right)^n&\ \mbox{if}\ x>0
\end{array}\right.\\
=\left\{\begin{array}{ll}
\beta\left(\frac{F(t)+F^\prime(t-)\frac{x}{n}+o(1/n)}{F(t)}\right)^n&\ \mbox{if}\ x<0\\
\beta&\ \mbox{if}\ x=0\\
1-(1-\beta)\left(\frac{1-G(t)-G^\prime(t+)\frac{x}{n}+o(1/n)}{1-G(t)}\right)^n&\ \mbox{if}\ x>0
\end{array}\right.\\
=\left\{\begin{array}{ll}
\beta\left(1+\frac{F^\prime(t-)}{F(t)}\frac{x}{n}+o(1/n)\right)^n&\ \mbox{if}\ x<0\\
\beta&\ \mbox{if}\ x=0\\
1-(1-\beta)\left(1-\frac{G^\prime(t+)}{1-G(t)}\frac{x}{n}+o(1/n)\right)^n&\ \mbox{if}\ x>0
\end{array}\right.\longrightarrow H(x)\quad\mbox{for every}\ x\in\mathbb{R},
\end{multline*}
where $H$ is defined by
$$H(x):=\left\{\begin{array}{ll}
\beta\exp\left(\frac{F^\prime(t-)}{F(t)}x\right)&\ \mbox{if}\ x\leq 0\\
1-(1-\beta)\exp\left(-\frac{G^\prime(t+)}{1-G(t)}x\right)&\ \mbox{if}\ x>0.
\end{array}\right.$$
So we have: an exponential distribution with mean $\frac{1-G(t)}{G^\prime(t+)}$ on $(0,\infty)$ with weight $1-\beta$; an exponential distribution with mean $\frac{F(t)}{F^\prime(t-)}$ on $(-\infty,0)$ with weight $\beta$.
\end{assertion}

Now we prove the LDP concerning the convergence to zero in $\mathbf{R1}$ of Assertion 
\ref{claim:ncMD-ex-DD-paper}.

\begin{proposition}\label{prop:LD-DD-paper}
Let $\{C_n:n\geq 1\}$ be the sequence defined by \eqref{eq:DD-paper-sequence}. Then 
$\{C_n:n\geq 1\}$ satisfies the LDP with speed $n$ and rate function 
$I_{\mathrm{LD}}$ defined by
$$I_{\mathrm{LD}}(x):=\left\{\begin{array}{ll}
\infty&\ \mbox{if}\ y\in (-\infty,-t]\\
-\log\frac{F(x+t)}{F(t)}&\ \mbox{if}\ x\in (-t,0]\\
-\log\frac{1-G(x+t)}{1-G(t)}&\ \mbox{if}\ x\in (0,\infty).
\end{array}\right.$$
\end{proposition}
\begin{proof}
We apply Lemma \ref{lemma}, with $s_n=n$, $I=I_{\mathrm{LD}}$ and, obviously, 
$\{C_n:n\geq 1\}$ as in \eqref{eq:DD-paper-sequence}.

\noindent\textsc{Proof of \eqref{eq:halfline1}.} For every $x>0$ we have
\begin{multline*}
\limsup_{n\to\infty}\frac{1}{n}\log P(C_n\geq x)
=\limsup_{n\to\infty}\frac{1}{n}\log(1-P(Z_n<x+t))\\
=\limsup_{n\to\infty}\frac{1}{n}\log\left((1-\beta)\left(\frac{1-G(x+t)}{1-G(t)}\right)^n\right)
=\log\left(\frac{1-G(x+t)}{1-G(t)}\right)=-I_{\mathrm{LD}}(x).
\end{multline*}

\noindent\textsc{Proof of \eqref{eq:halfline2}.} For every $x<0$ we have (when $x\leq-t$ we use the convention
$\log 0=-\infty$)
\begin{multline*}
\limsup_{n\to\infty}\frac{1}{n}\log P(C_n\leq x)
=\limsup_{n\to\infty}\frac{1}{n}\log P(Z_n\leq x+t)\\
=\limsup_{n\to\infty}\frac{1}{n}\log\left(\beta\left(\frac{F(x+t)}{F(t)}\right)^n\right)
=\log\left(\frac{F(x+t)}{F(t)}\right)=-I_{\mathrm{LD}}(x).
\end{multline*}

\noindent\textsc{Proof of \eqref{eq:local-bound}.} We want to show that, for every $x>-t$ and for every open 
set $O$ such that $x\in O$, we have
$$\liminf_{n\to\infty}\frac{1}{n}\log P(C_n\in O)\geq\left\{\begin{array}{ll}
\log\frac{F(x+t)}{F(t)}&\ \mbox{if}\ x\in (-t,0]\\
\log\frac{1-G(x+t)}{1-G(t)}&\ \mbox{if}\ x\in (0,\infty).
\end{array}\right.$$
The case $x=0$ is immediate; indeed, since $C_n$ converges in probability to zero, we have trivially 
$0\geq 0$. So, from now on, we suppose that $x>-t$ with $x\neq 0$ and we have two cases: $x\in(-t,0)$
and $x\in(0,\infty)$.

If $x\in(-t,0)$ we take $\delta>0$ small enough to have $(x-\delta,x+\delta)\subset O\cap(-t,0)$; then
\begin{multline*}
P(C_n\in O)\geq P(C_n\in(x-\delta,x+\delta))=P(Z_n<x+\delta+t)-P(Z_n\leq x-\delta+t)\\
=\beta\left(\left(\frac{F(x+\delta+t)}{F(t)}\right)^n-\left(\frac{F(x-\delta+t)}{F(t)}\right)^n\right)\\
=\frac{\beta}{F(t)}\underbrace{(F(x+\delta+t)-F(x-\delta+t))}_{>0}
\underbrace{\sum_{k=0}^{n-1}\left(\frac{F(x+\delta+t)}{F(t)}\right)^k
\left(\frac{F(x-\delta+t)}{F(t)}\right)^{n-1-k}}_{\geq\left(\frac{F(x+\delta+t)}{F(t)}\right)^{n-1}>0},
\end{multline*}
and therefore
$$\liminf_{n\to\infty}\frac{1}{n}\log P(C_n\in O)\geq\log\left(\frac{F(x+\delta+t)}{F(t)}\right).$$
Then we get \eqref{eq:local-bound} for $x\in(-t,0)$ by letting $\delta$ go to zero (by the continuity
of $F$).

If $x\in(0,\infty)$ we take $\delta>0$ small enough to have $(x-\delta,x+\delta)\subset O\cap(0,\infty)$; then
\begin{multline*}
P(C_n\in O)\geq P(C_n\in(x-\delta,x+\delta))=P(Z_n<x+\delta+t)-P(Z_n\leq x-\delta+t)\\
=1-(1-\beta)\left(\frac{1-G(x+\delta+t)}{1-G(t)}\right)^n-\left(1-(1-\beta)\left(\frac{1-G(x-\delta+t)}{1-G(t)}\right)^n\right)\\
=(1-\beta)\left(\left(\frac{1-G(x-\delta+t)}{1-G(t)}\right)^n-\left(\frac{1-G(x+\delta+t)}{1-G(t)}\right)^n\right)\\
=\underbrace{\frac{1-\beta}{1-G(t)}(G(x+\delta+t)-G(x-\delta+t))}_{>0}\underbrace{\sum_{k=0}^{n-1}
\left(\frac{1-G(x-\delta+t)}{1-G(t)}\right)^k\left(\frac{1-G(x+\delta+t)}{1-G(t)}\right)^{n-1-k}}
_{\geq\left(\frac{1-G(x-\delta+t)}{1-G(t)}\right)^{n-1}>0}
\end{multline*}
and therefore
$$\liminf_{n\to\infty}\frac{1}{n}\log P(C_n\in O)\geq\log\left(\frac{1-G(x-\delta+t)}{1-G(t)}\right).$$
Then we get \eqref{eq:local-bound} for $x\in(0,\infty)$ by letting $\delta$ go to zero (by the continuity
of $G$).
\end{proof}

Now we prove the moderate deviation result. We shall see that, for this example, the family of positive 
numbers $\{a_n:n\geq 1\}$ such that \eqref{eq:MD-conditions} holds (i.e. $a_n\to 0$ and $a_nn\to\infty$) 
has to verify also the stricter condition
\begin{equation}\label{eq:restriction}
a_n\log n\to 0.
\end{equation}
Obviously \eqref{eq:restriction} yields $a_n\to 0$.

\begin{proposition}\label{prop:MD-DD-paper}
For every family of positive numbers $\{a_n:n\geq 1\}$ such that \eqref{eq:MD-conditions} and 
\eqref{eq:restriction} hold (i.e. $a_n\log n\to 0$ and $a_nn\to\infty$), the sequence of random variables 
$\{a_nn(Z_n-t):n\geq 1\}$ satisfies the LDP with speed $1/a_n$ and rate function $I_{\mathrm{MD}}$ defined 
by
$$I_{\mathrm{MD}}(x):=\left\{\begin{array}{ll}
-\frac{F^\prime(t-)}{F(t)}x&\ \mbox{if}\ x\leq 0\\
\frac{G^\prime(t+)}{1-G(t)}x&\ \mbox{if}\ x>0.
\end{array}\right.$$
\end{proposition}
\begin{proof}
We apply Lemma \ref{lemma} for every choice of the positive scalings $\{a_n:n\geq 1\}$ as
in the statement of the proposition, with $s_n=1/a_n$, $I=I_{\mathrm{MD}}$ and 
$\{C_n:n\geq 1\}$ as in \eqref{eq:DD-paper-sequence}.

\noindent\textsc{Proof of \eqref{eq:halfline1}.} For every $x>0$ we have
\begin{multline*}
\limsup_{n\to\infty}a_n\log P(a_nn(Z_n-t)\geq x)
=\limsup_{n\to\infty}a_n\log\left(1-P\left(Z_n<t+\frac{x}{a_nn}\right)\right)\\
=\limsup_{n\to\infty}a_n\log\left((1-\beta)\left(\frac{1-G\left(t+\frac{x}{a_nn}\right)}{1-G(t)}\right)^n\right)
=\limsup_{n\to\infty}a_nn\log\left(\frac{1-G\left(t+\frac{x}{a_nn}\right)}{1-G(t)}\right)\\
=\limsup_{n\to\infty}a_nn\log\left(\frac{1-G(t)-G^\prime(t+)\frac{x}{a_nn}+o(\frac{1}{a_nn})}{1-G(t)}\right)\\
=\limsup_{n\to\infty}a_nn\log\left(1-\frac{G^\prime(t+)}{1-G(t)}\frac{x}{a_nn}+o\left(\frac{1}{a_nn}\right)\right)
=-\frac{G^\prime(t+)}{1-G(t)}x=-I_{\mathrm{MD}}(x).
\end{multline*}
\noindent\textsc{Proof of \eqref{eq:halfline2}.} For every $x<0$ we have (note that $t+\frac{x}{a_nn}>0$ for $n$ 
large enough)
\begin{multline*}
\limsup_{n\to\infty}a_n\log P(a_nn(Z_n-t)\leq x)
=\limsup_{n\to\infty}a_n\log P\left(Z_n\leq t+\frac{x}{a_nn}\right)\\
=\limsup_{n\to\infty}a_n\log\left(\beta\left(\frac{F\left(t+\frac{x}{a_nn}\right)}{F(t)}\right)^n\right)
=\limsup_{n\to\infty}a_nn\log\left(\frac{F\left(t+\frac{x}{a_nn}\right)}{F(t)}\right)\\
=\limsup_{n\to\infty}a_nn\log\left(\frac{F(t)+F^\prime(t-)\frac{x}{a_nn}+o(1/(a_nn))}{F(t)}\right)\\
=\limsup_{n\to\infty}a_nn\log\left(1+\frac{F^\prime(t-)}{F(t)}\frac{x}{a_nn}+o\left(\frac{1}{a_nn}\right)\right)
=\frac{F^\prime(t-)}{F(t)}x=-I_{\mathrm{MD}}(x).
\end{multline*}

\noindent\textsc{Proof of \eqref{eq:local-bound}.} We show that, for every $x\in\mathbb{R}$ and 
for every open set $O$ such that $x\in O$, we have
$$\liminf_{n\to\infty}a_n\log P(a_nn(Z_n-t)\in O)\geq\left\{\begin{array}{ll}
\frac{F^\prime(t-)}{F(t)}x&\ \mbox{if}\ x\leq 0\\
-\frac{G^\prime(t+)}{1-G(t)}x&\ \mbox{if}\ x>0.
\end{array}\right.$$
The case $x=0$ is immediate; indeed, since $a_nn(Z_n-t)$ converges in probability to
zero by the Slutsky Theorem (by $a_n\to 0$ and the weak convergence in $\mathbf{R2}$ in Assertion 
\ref{claim:ncMD-ex-DD-paper}), we have 
trivially $0\geq 0$. So, from now on, we suppose that $x\neq 0$ and we have two 
cases: $x\in(-\infty,0)$ and $x\in(0,\infty)$. 

If $x\in(-\infty,0)$ we take $\delta>0$ small enough to have $(x-\delta,x+\delta)\subset O\cap(-\infty,0)$; then
(note that $t+\frac{x\pm\delta}{a_nn}>0$ for $n$ large enough)
\begin{multline*}
P(a_nn(Z_n-t)\in O)\geq P(a_nn(Z_n-t)\in(x-\delta,x+\delta))
=P\left(Z_n<t+\frac{x+\delta}{a_nn}\right)-P\left(Z_n\leq t+\frac{x-\delta}{a_nn}\right)\\
=\frac{\beta}{(F(t))^n}\left(\left(F\left(t+\frac{x+\delta}{a_nn}\right)\right)^n-
\left(F\left(t+\frac{x-\delta}{a_nn}\right)\right)^n\right)\\
=\frac{\beta}{(F(t))^n}\underbrace{\left(F\left(t+\frac{x+\delta}{a_nn}\right)-
F\left(t+\frac{x-\delta}{a_nn}\right)\right)}_{=F^\prime(t-)\frac{2\delta}{a_nn}+o\left(\frac{1}{a_nn}\right)}\\
\cdot\underbrace{\sum_{k=0}^{n-1}\left(F\left(t+\frac{x+\delta}{a_nn}\right)\right)^k
\left(F\left(t+\frac{x-\delta}{a_nn}\right)\right)^{n-1-k}}_{\geq\left(F\left(t+\frac{x+\delta}{a_nn}\right)\right)^{n-1}},
\end{multline*}
and therefore (in the final step we use the condition $a_n\log n\to 0$ in \eqref{eq:restriction})
\begin{multline*}
\liminf_{n\to\infty}a_n\log P(a_nn(Z_n-t)\in O)\geq 
\liminf_{n\to\infty}a_n\log\left\{\left(F^\prime(t-)\frac{2\delta}{a_nn}+o\left(\frac{1}{a_nn}\right)\right)
\left(\frac{F\left(t+\frac{x+\delta}{a_nn}\right)}{F(t)}\right)^{n-1}\right\}\\
=\liminf_{n\to\infty}\left\{a_n\log\left(F^\prime(t-)\frac{2\delta}{a_nn}+o\left(\frac{1}{a_nn}\right)\right)
+a_n(n-1)\log\left(\frac{F\left(t+\frac{x+\delta}{a_nn}\right)}{F(t)}\right)\right\}\\
=\liminf_{n\to\infty}\left\{-a_n\log n+a_n(n-1)\log\left(\frac{F(t)+F^\prime(t-)\frac{x+\delta}{a_nn}+o(1/(a_nn))}{F(t)}\right)\right\}\\
=\liminf_{n\to\infty}\left\{-a_n\log n+a_n(n-1)\log\left(1+\frac{F^\prime(t-)}{F(t)}\frac{x+\delta}{a_nn}+o\left(\frac{1}{a_nn}\right)\right)\right\}
=\frac{F^\prime(t-)}{F(t)}(x+\delta).
\end{multline*}
So we obtain \eqref{eq:local-bound} for $x\in(-\infty,0)$ by letting $\delta$ go to zero.

If $x\in(0,\infty)$ we take $\delta>0$ small enough to have $(x-\delta,x+\delta)\subset O\cap(0,\infty)$; then
\begin{multline*}
P(a_nn(Z_n-t)\in O)\geq P(a_nn(Z_n-t)\in(x-\delta,x+\delta))
=P\left(Z_n<t+\frac{x+\delta}{a_nn}\right)-P\left(Z_n\leq t+\frac{x-\delta}{a_nn}\right)\\
=\frac{1-\beta}{(1-G(t))^n}\left(\left(1-G\left(t+\frac{x-\delta}{a_nn}\right)\right)^n
-\left(1-G\left(t+\frac{x+\delta}{a_nn}\right)\right)^n\right)\\
=\frac{1-\beta}{(1-G(t))^n}\underbrace{\left(G\left(t+\frac{x+\delta}{a_nn}\right)-
G\left(t+\frac{x-\delta}{a_nn}\right)\right)}_{=G^\prime(t+)\frac{2\delta}{a_nn}+o\left(\frac{1}{a_nn}\right)}\\
\cdot\underbrace{\sum_{k=0}^{n-1}\left(1-G\left(t+\frac{x-\delta}{a_nn}\right)\right)^k
\left(1-G\left(t+\frac{x+\delta}{a_nn}\right)\right)^{n-1-k}}_{\geq\left(1-G\left(t+\frac{x-\delta}{a_nn}\right)\right)^{n-1}},
\end{multline*}
and therefore (in the final step we use the condition $a_n\log n\to 0$ in \eqref{eq:restriction})
\begin{multline*}
\liminf_{n\to\infty}a_n\log P(a_nn(Z_n-t)\in O)\geq 
\liminf_{n\to\infty}a_n\log\left\{\left(G^\prime(t+)\frac{2\delta}{a_nn}+o\left(\frac{1}{a_nn}\right)\right)
\left(\frac{1-G\left(t+\frac{x-\delta}{a_nn}\right)}{1-G(t)}\right)^{n-1}\right\}\\
=\liminf_{n\to\infty}\left\{a_n\log\left(G^\prime(t+)\frac{2\delta}{a_nn}+o\left(\frac{1}{a_nn}\right)\right)
+a_n(n-1)\log\left(\frac{1-G\left(t+\frac{x-\delta}{a_nn}\right)}{1-G(t)}\right)\right\}\\
=\liminf_{n\to\infty}\left\{-a_n\log n+a_n(n-1)\log\left(\frac{1-G(t)-G^\prime(t+)\frac{x-\delta}{a_nn}+o(\frac{1}{a_nn})}{1-G(t)}\right)\right\}\\
=\liminf_{n\to\infty}\left\{-a_n\log n+a_n(n-1)\log\left(1-\frac{G^\prime(t+)}{1-G(t)}\frac{x-\delta}{a_nn}
+o\left(\frac{1}{a_nn}\right)\right)\right\}
=-\frac{G^\prime(t+)}{1-G(t)}(x-\delta).
\end{multline*}
So we obtain \eqref{eq:local-bound} for $x\in(0,\infty)$ by letting $\delta$ go to zero.
\end{proof}

\subsection*{Funding}
This work has been supported by MIUR Excellence Department Project awarded to the Department of Mathematics,
University of Rome Tor Vergata (CUP E83C18000100006), by University of Rome Tor Vergata (project "Asymptotic 
Methods in Probability" (CUP E89C20000680005) and project "Asymptotic Properties in Probability" 
(CUP E83C22001780005)) and by Indam-GNAMPA.

\end{document}